\documentclass[10pt,fullpage,leqno]{article}

\usepackage{amsfonts}
\usepackage{amsmath}
\usepackage{epsfig}
\usepackage{pstricks}
\usepackage{graphicx}
\usepackage{enumerate}

\def\R{\hbox{\bf R}}
\def\Z{\hbox{\bf Z}}

\def\N{\hbox{\bf N}}

\newcommand{\ba}{\begin{eqnarray}}
\newcommand{\ea}{\end{eqnarray}}

\newtheorem{theo}{\bf Theorem}[section]
\newtheorem{lem}[theo]{\bf Lemma}
\newtheorem{pro}[theo]{\bf Proposition}
\newtheorem{cor}[theo]{\bf Corollary}
\newtheorem{defi}[theo]{\bf Definition}
\newtheorem{rem}[theo]{\bf Remark}

\renewcommand{\N}{{\mathbb N}}
\renewcommand{\R}{{\mathbb R}}

\renewcommand{\Z}{{\mathbb Z}}

\newenvironment{Proofc}[1]{\smallskip\par\noindent\textsc{#1} }%
  {\hfill$\Box$\bigskip\par}
\newenvironment{Proof}{\begin{Proofc}{\textbf{Proof}}}{\end{Proofc}}

\newcommand{\argmax}[1]{ \smash{\mathop{{\rm argmax}}\limits_{#1}}\ }
\DeclareMathOperator*{\esssup}{ess\,sup}

\nonstopmode

\begin{document}

\title{\bf A convergent scheme for Hamilton-Jacobi equations\\ on a junction: application to traffic}

\author{G. Costeseque\footnotemark[1]  \footnotemark[2] , 
J-P. Lebacque\footnotemark[2] , 
R. Monneau\footnotemark[1]}

\footnotetext[1]{Universit\'e Paris-Est, Ecole des Ponts ParisTech, CERMICS, 6 \& 8 avenue Blaise
  Pascal, Cit\'e Descartes, Champs sur Marne, 77455 Marne la Vall\'ee
  Cedex 2, France.}

\footnotetext[2]{Universit\'e Paris-Est, IFSTTAR, GRETTIA, 14-20 Boulevard Newton,
Cit\'e Descartes, Champs sur Marne, 77447 Marne la Vall\'ee Cedex 2, France.}
  

\maketitle

\vspace{15pt}




 \centerline{\small{\bf{Abstract}}}
{\small{In this paper, we consider first order Hamilton-Jacobi (HJ) equations posed on a ``junction'', that is to say the union of a finite number of half-lines with a unique common point.
For this continuous HJ problem, we propose a finite difference scheme and prove two main results. 
As a first result, we show bounds on the discrete gradient and time derivative of the numerical solution.
Our second result is the convergence (for a subsequence) of the numerical solution towards a viscosity solution of the continuous HJ problem, as the mesh size goes to zero.
When the solution of the continuous HJ problem is unique, we recover the full convergence of the numerical solution.
We apply this scheme to compute the densities of cars for a traffic model. We recover the well-known Godunov scheme outside the junction point and we give a numerical illustration.}}\hfill\break

 \noindent{\small{\bf{Keywords:}}} {\small{Hamilton-Jacobi equations, junctions, viscosity solutions, numerical scheme, traffic problems.}}\hfill\break
 \noindent{\small{\bf{MSC Classification:}}} {\small{65M12, 65M06, 35F21, 90B20.}}\hfill\break

\parindent 0cm


\section{Introduction}
\label{s1}

The main goal of this paper is to prove properties of a numerical scheme to solve Hamilton-Jacobi (HJ) equations posed on a junction.
We also propose a traffic application that can be directly found in Section \ref{sect::application}.

\subsection{Setting of the PDE problem}

In this subsection, we first define the junction, then the space of functions on the junction and finally the Hamilton-Jacobi equations.
We follow \cite{IMZ}.\\

\begin{figure}[!ht]
\begin{center}
\resizebox{3cm}{!}{\input{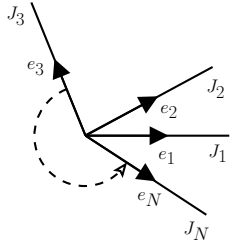}}
\caption{Junction model}
\label{Junction model}
\end{center}
\end{figure}

\textbf{The junction.} Let us consider $N \geq 1$ different unit vectors $e_{\alpha} \in \R^2$ for ${\alpha}=1,\ldots,N$. We define the branches as the half-lines generated by these unit vectors
$$ J_{\alpha}=[0,+\infty) e_{\alpha} \quad \text{and} \quad J^*_{\alpha}=J_{\alpha} \setminus \{0_{\R^2}\}, \quad \text{for all} \quad {\alpha}=1,\ldots,N, $$
and the whole \textit{junction} (see Figure \ref{Junction model}) as
$$ J=\bigcup_{{\alpha}=1,\ldots,N} J_{\alpha} .$$
The origin $y=0_{\R^2}$ (we just call it ``$y=0$'' in the following) is called the \textit{junction point}. For a time $T>0$, we also consider the time-space domain defined as
$$ J_T=(0,T)\times J .$$

\textbf{Space of test functions.}  For a function $u : J_T \rightarrow \R $, we denote by $u^{\alpha}$ the ``restriction'' of $u$ to $ (0,T) \times J_{\alpha}$ defined as follows for $x\ge 0$
$$u^{\alpha}(t,x):=u(t,x e_{\alpha}).$$
Then we define the natural space of functions on the junction:
\begin{equation}\label{spacefunct}
C^1_*(J_T)= \left \{ u \in C(J_T), \quad u^{\alpha} \in C^1 \left( (0,T) \times [0,+\infty) \right) \quad \mbox{for} \quad {\alpha}=1,\ldots,N \right \}.
\end{equation}

In particular for $u \in C^1_*(J_T)$ and $y=x e_{\alpha}$ with $x \geq 0$, we define
$$ u_t(t,y)=u^{\alpha}_t(t,x)=\dfrac{\partial u^{\alpha}}{\partial t}(t,x) \quad \mbox{and} \quad  u^{\alpha}_x(t,x)=\dfrac{\partial u^{\alpha}}{\partial x}(t,x) .$$

\textbf{HJ equation on the junction.}  We are interested in continuous functions $u: [0,T) \times J \rightarrow \R$ which are \textit{viscosity solutions} (see Definition~\ref{defviscosity}) on $J_T$ of
\begin{equation}\label{HJ}
\begin{cases}
u^{\alpha}_t + H_{\alpha}(u^{\alpha}_x) =0  \quad &\mbox{on} \quad (0,T)\times (0,+\infty),  \quad \mbox{for} \quad {\alpha}=1,\ldots ,N , \\
\\
\left.\begin{array}{l}
u^{\beta} =: u, \quad \text{for all} \quad \beta = 1,\ldots,N \\
\\
u_t + \displaystyle{\max_{\beta=1,\ldots ,N}} \ H^{-}_{\beta}(u^{\beta}_x) =0
\end{array}\right| \quad &\mbox{on} \quad (0,T)\times \{ 0 \} ,
\end{cases}
\end{equation}
for functions $H_\alpha$ and $H_\alpha^-$ that will be defined below in assumption (A1).

We consider an initial condition
\begin{equation}\label{CI}
u^{\alpha}(0,x)=u^{\alpha}_0(x), \quad \mbox{with} \quad x \in [0,+\infty) \quad \mbox{for} \quad {\alpha}=1,\ldots ,N.
\end{equation}

We make the following assumptions:\\
\noindent (A0) {\bf Initial data}\\
The initial data $u_0 :=\left( u_0^\alpha \right)_{\alpha}$ is globally Lipschitz continuous on $J$, i.e. each associated $u^{\alpha}_0$ is Lipschitz continuous on $[0,+\infty)$ and $u^{\alpha}_0(0)=u^{\beta}_0(0)$ for any $\alpha\not= \beta$.\\

\noindent (A1) {\bf Hamiltonians}\\
For each $\alpha=1,...,N$,
\begin{list}{$\bullet$}{}
\item we consider functions $H_\alpha \in C^1(\R;\R)$ which are coercive, i.e. $ \displaystyle{\lim_{|p|\to +\infty}} H_{\alpha}(p) = + \infty$;
\item we assume that there exists a  $p^{\alpha}_0 \in \R$ such that $H_{\alpha}$ is non-increasing on $(- \infty,p^{\alpha}_0]$ and non-decreasing on $[p^{\alpha}_0, +\infty)$, and we set:
\begin{equation}\label{Hamiltonians+-}
H^{-}_{\alpha}(p) = \begin{cases} 
H_{\alpha}(p) \quad &\mbox{for} \quad p \leq p_{0}^ {\alpha} \\
\\
H_{\alpha}(p_{0}^ {\alpha}) \quad &\mbox{for} \quad p \geq p_{0}^ {\alpha} \quad
\end{cases}
\quad \mbox{and}\quad 
H^{+}_{\alpha}(p) = \begin{cases} 
H_{\alpha}(p_{0}^ {\alpha}) \quad &\mbox{for} \quad p \leq p_{0}^ {\alpha} \\
\\
H_{\alpha}(p) \quad &\mbox{for} \quad p \geq p_{0}^ {\alpha} \quad
\end{cases}
\end{equation}
where $H_\alpha^-$ is non-increasing and $H_\alpha^+$ is non-decreasing.
\end{list}

\begin{rem}
In assumption (A1), we assume that $p^\alpha_0$ is unique, i.e. there is no plateau at the minimum of $H_\alpha$.
This condition is not fundamental, but simplifies the presentation of the work.
\end{rem}

\subsection{Presentation of the scheme}

We denote by $\Delta x$ the space step and by $\Delta t$ the time step. We denote by $U^{\alpha,n}_i$ an approximation of $ u^{\alpha}(n \Delta t, i \Delta x)$ for $n \in \N, \ i \in \N$, where $\alpha$ stands for the index of the considered branch. 

We define the discrete space derivatives
\begin{equation}\label{defp}
p^{\alpha,n}_{i,+}:=\dfrac{U^{\alpha,n}_{i+1} - U^{\alpha,n}_i}{\Delta x} \quad \mbox{and} \quad p^{\alpha,n}_{i,-}:=\dfrac{U^{\alpha,n}_i - U^{\alpha,n}_{i-1}}{\Delta x} ,
\end{equation}
and similarly the discrete time derivative
\begin{equation}\label{defW}
W^{\alpha,n}_{i}:=\dfrac{U^{\alpha,n+1}_{i} - U^{\alpha,n}_i}{\Delta t} .
\end{equation}

Then we consider the following numerical scheme corresponding to the discretization of the HJ equation (\ref{HJ}) for $n\ge 0$:
\begin{equation}\label{NS}
\begin{cases}
\dfrac{U^{\alpha ,n+1}_{i} -U^{\alpha ,n}_{i}}{\Delta t} + \max \left\{ H^{+}_{\alpha} (p^{\alpha ,n}_{i,-}), H^{-}_{\alpha} (p^{\alpha ,n}_{i,+}) \right\} =0,& \quad \mbox{for}\quad  i \ge 1, \quad \alpha = 1,\ldots,N , \\
\\
\left.\begin{array}{l}
U^{\beta,n}_{0} =: U^{n}_{0}, \quad \text{for all} \quad \beta = 1,\ldots,N \\
\\
\dfrac{U^{n+1}_{0} -U^{n}_{0}}{\Delta t} + \displaystyle{\max_{\beta = 1,\ldots ,N}} \ H^{-}_{\beta} (p^{\beta ,n}_{0,+}) =0
\end{array}\right|& \quad \mbox{for} \quad i=0,
\end{cases}
\end{equation}
with the initial condition
\begin{equation}\label{CINS}
U^{\alpha,0}_{i} = u^{\alpha}_{0} (i \Delta x) \quad \mbox{for} \quad i \ge 0, \quad \alpha = 1,\ldots,N.
\end{equation}

It is natural to introduce the following Courant-Friedrichs-Lewy (CFL) condition \cite{CFL}:
\begin{equation}\label{CFL}
\dfrac{\Delta x}{\Delta t} \geq \sup_{\substack{ \alpha=1,\dots ,N \\ i \geq 0, \ 0 \leq n \leq n_T}} | H'_{\alpha}(p^{\alpha,n}_{i,+}) |
\end{equation}
where the integer $n_T$ is assumed to be defined as $ n_T= \left \lfloor \dfrac{T}{\Delta t} \right \rfloor $ for a given $T >0$.\\

We then have

\begin{pro}\label{pro1} {\bf (Monotonicity of the numerical scheme)}\\
Let $U^n:=\left( U^{\alpha,n}_i \right)_{\alpha,i}$ and $V^n:=\left( V^{\alpha,n}_i \right)_{\alpha,i}$ two solutions of (\ref{NS}).
If the CFL condition (\ref{CFL}) is satisfied and if $U^0 \leq V^{0}$, then the numerical scheme (\ref{NS}) is monotone, that is
$$U^n \leq V^n \quad \mbox{for any} \quad n \in \{0,...,n_T \} .$$
\end{pro}

Our scheme (\ref{NS}) is related to the Godunov scheme for conservation laws in one space dimension,
as it is explained in our application to traffic in Section \ref{sect::application}.

\subsection{Main results}
\label{mainresults}

We first notice that it is not obvious to satisfy  the CFL condition (\ref{CFL}) because for any $\alpha$, $i$ and $n$, the discrete gradients $p^{\alpha ,n}_{i,+} $ depends itself on $ \Delta t $ through the scheme (\ref{NS}). For this reason, we will consider below a more restrictive CFL condition (see (\ref{CFL'})) that can be checked from the initial data. 
To this end, we need to introduce a few notations.

For sake of clarity we first consider $\sigma \in \left\{+1,-1\right\}$ denoted by abuse of notation $\sigma \in \left\{+,-\right\}$ in the remaining, with the convention $-\sigma=-$ if $\sigma=+$ and $-\sigma = +$ if $\sigma=-$.

Under assumption (A1), we need to use a sort of inverse of  $(H_\alpha^\pm)$ that we define naturally for 
$\sigma \in \left\{+,-\right\}$ as:
\begin{equation}\label{rem:1}
(H^{-\sigma}_{\alpha})^{-1} (a) := \sigma \left(\inf \{\sigma p,\ \ H^{-\sigma}_{\alpha} (p)=a \}\right)
\end{equation}
with the additional convention that $(H^{\pm}_{\alpha})^{-1} (+\infty)=\pm\infty$.

We set 
\begin{equation}\label{defcontgrad}
\begin{cases}
\underline{p}_{\alpha} =  (H^{-}_{\alpha})^{-1} (-{m}^{0})\\
\\
\overline{p}_{\alpha} =  (H^{+}_{\alpha})^{-1} (-{m}^{0})
\end{cases}
\quad \mbox{with} \quad 
{m}^{0} = \displaystyle{\inf_{\substack{\beta =1,...,N ,\\
i\in\N}}} \ W^{\beta,0}_i
\end{equation}
where $(W^{\beta,0}_i)_{\beta,i}$, defined in (\ref{defW}), is given by the scheme (\ref{NS}) for $n=0$ in terms of 
$(U^{\beta,0}_i)_{\beta,i}$ (itself defined in (\ref{CINS})). 
It is important to notice that with this construction,  $\underline{p}_{\alpha}$ and $\overline{p}_{\alpha}$ depend on $\Delta x$, but not on $\Delta t$.

We now consider a more restrictive CFL condition given by
\begin{equation}\label{CFL'}
\dfrac{\Delta x}{\Delta t} \geq \sup_{\substack{\alpha=1,\dots ,N \\ p_{\alpha} \in [\underline{p}_{\alpha}, \overline{p}_{\alpha}]}} | H'_{\alpha}(p_{\alpha}) |
\end{equation}
which is then satisfied for $\Delta t$ small enough.

Our first main result is the following:

\begin{theo}\label{gradientestimate} {\bf (Gradient and time derivative estimates)}\\
Assume (A1). 
If $(U^{\alpha,n}_i)$ is the numerical solution of (\ref{NS})-(\ref{CINS}) and if the CFL condition (\ref{CFL'}) is satisfied with $m^0$ finite, then  the following two properties hold for any $n \geq 0$:
\begin{enumerate}[(i)]
\item For $\underline{p}_{\alpha}$ and $\overline{p}_{\alpha}$ defined in (\ref{defcontgrad}), we have the following gradient estimate:
\begin{equation}\label{eq::g1}
\underline{p}_{\alpha} \leq  p^{\alpha,n}_{i,+} \leq  \overline{p}_{\alpha}, \quad \text{for all} \quad i \geq 0, \quad \text{and} \quad \alpha=1,...,N.
\end{equation}
\item Considering ${M}^n= \displaystyle{\sup_{\alpha,i}} \ W^{\alpha,n}_i$ and ${m}^n= \displaystyle{\inf_{\alpha,i}} \ W^{\alpha,n}_i$, we have the following time derivative estimate:
\begin{equation}\label{eq::g2}
{m}^{0} \leq  {m}^{n} \leq  {m}^{n+1} \leq  {M}^{n+1}  \leq  {M}^{n}  \leq  {M}^{0}.
\end{equation}
\end{enumerate}
\end{theo}

\begin{rem}
\label{rem:2}
Notice that due to (\ref{eq::g1}), the more restrictive CFL condition (\ref{CFL'}) implies the natural CFL condition (\ref{CFL}) for any $n_T \geq 0$.
\end{rem}

Our second main result is the following convergence result which also gives the existence of a solution to equations (\ref{HJ})-(\ref{CI}).

\begin{theo}\label{convergence1} {\bf (Convergence of the numerical solution up to a subsequence)}\\
Assume (A0)-(A1). Let $T>0$ and $\varepsilon = (\Delta t, \Delta x)$ such that the CFL condition (\ref{CFL'}) is satisfied. If $u :=\left( u^{\alpha} \right)_{\alpha}$ is a solution of (\ref{HJ})-(\ref{CI}) in the sense of Definition~\ref{defviscosity}, then there exist a subsequence $\varepsilon'$ of $\varepsilon$ such that the numerical solution $(U^{\alpha,n}_i)$ of (\ref{NS})-(\ref{CINS}) converges to $u$ when $\varepsilon$ goes to zero, locally uniformly on any compact set $\mathcal{K} \subset [0,T) \times J$, i.e.
\begin{equation}\label{eq::limit1}
\limsup_{\varepsilon\to 0} \ \sup_{(n \Delta t,i \Delta x) \in \mathcal{K}} \ | u^{\alpha} (n \Delta t,i \Delta x) - U^{\alpha, n}_i | = 0,
\end{equation}
where the index $\alpha$ in (\ref{eq::limit1}) is chosen such that $(n \Delta t,i \Delta x) \in \mathcal{K} \cap [0,T) \times J_{\alpha}$.
\end{theo}

In order to give below sharp Lipschitz estimates on the continuous solution $u$, we first define $L^{\alpha,-}$ and $L^{\alpha,+}$ as 
the best Lipschitz constants for the initial data $u^{\alpha}_0$, i.e. satisfying for any $x\geq 0$ and $a \geq 0$ 
\begin{equation}\label{LipConst}
a L^{\alpha,-} \leq u^{\alpha}_0 (x+a) - u^{\alpha}_0 (x) \leq a L^{\alpha,+}.
\end{equation}

Let us consider
\begin{equation}\label{defcontgrad1}
\begin{cases}
m^0_0:= \displaystyle{\inf_{\substack{\alpha = 1,...,N \\ L^{\alpha,-} \leq p_{\alpha} \leq L^{\alpha,+}}}} \ -H_{\alpha}(p_{\alpha}), \\
M^0_0:= \displaystyle{\max \left[ \max_{\alpha=1,...,N} \ \left\{-\max_{\sigma \in \{ +,- \}} \ H^{-\sigma}_{\alpha}(L^{\alpha,\sigma}) \right\} , \ -\max_{\alpha=1,...,N} \ \bigg\{ H^{-}_{\alpha}(L^{\alpha,+}) \bigg\} \right] },
\end{cases}
\end{equation}
and
\begin{equation}\label{defcontgrad2}
\begin{cases}
\underline{p}_{\alpha}^0:= (H^{-}_{\alpha})^{-1}(-m^0_0), \\
\\
\overline{p}_{\alpha}^0:= (H^{+}_{\alpha})^{-1}(-m^0_0).
\end{cases}
\end{equation}

\begin{cor}\label{cor2} {\bf (Gradient and time derivative estimates for a continuous solution)}\\
Assume (A0)-(A1). Let $T>0$ and $u:=\left( u^{\alpha} \right)_{\alpha}$ be a solution of (\ref{HJ})-(\ref{CI}) constructed in Theorem~\ref{convergence1}. 
Then for all $a \geq 0$, for all $0 \leq t \leq T$ and $x \geq 0$, the function $u$ satisfies the following properties:
\begin{equation}\label{contestimates}
\begin{cases}
a m^0_0 \leq u^{\alpha}(t+a,x)-u^{\alpha}(t,x) \leq a M^0_0, \\
\\
a \underline{p}_{\alpha}^0 \leq u^{\alpha}(t,x+a)-u^{\alpha}(t,x) \leq a \overline{p}_{\alpha}^0,
\end{cases}
\end{equation}
where $m^0_0$, $M^0_0$, $\underline{p}_{\alpha}^0$ and $\overline{p}_{\alpha}^0$ are defined in (\ref{defcontgrad1}) and (\ref{defcontgrad2}).
\end{cor}

Recall that under the general assumptions of Theorem~\ref{convergence1}, i.e. (A0)-(A1), 
the uniqueness of a solution $u$ of (\ref{HJ})-(\ref{CI}) is not known. 
If we replace condition (A1) by a stronger assumption (A1') below,
it is possible to recover the uniqueness of the solution (see \cite{IMZ} and Theorem~\ref{existence&uniqueness} below).\\
This is the following assumption:\\

\noindent (A1') {\bf Strong convexity}\\
There exists a constant $\gamma>0$, such that for each $\alpha=1,...,N$, 
there exists a lagrangian function $L_\alpha\in C^2(\R;\R)$ satisfying $L_\alpha'' \ge \gamma >0$ 
such that $H_\alpha$ is the Legendre-Fenchel transform of $L_\alpha$, i.e.
\begin{equation}\label{LF1}
H_{\alpha}(p)=L^*_{\alpha} (p) = \sup_{q \in \R} \ (pq - L_{\alpha}(q))
\end{equation}
and
\begin{equation}\label{LF2}
H^{-}_{\alpha} (p)= \sup_{q \leq 0} \ (pq - L_{\alpha}(q))  \quad \mbox{and} \quad H^{+}_{\alpha} (p)= \sup_{q \geq 0} \ (pq - L_{\alpha}(q)).
\end{equation}

We can easily check that assumption (A1') implies assumption (A1).

We are now ready to recall the following result extracted from \cite{IMZ}:

\begin{theo}\label{existence&uniqueness} {\bf (Existence and uniqueness for a solution of the HJ problem)}\\
Assume (A0)-(A1') and let $T>0$. Then there exists a unique viscosity solution $u$ of (\ref{HJ})-(\ref{CI}) on $J_T$ in the sense of the Definition~\ref{defviscosity}, satisfying for some constant $C_T>0$
$$ | u(t,y) - u_0(y)| \leq C_T \quad \mbox{for all} \quad (t,y) \in J_T. $$
Moreover the function $u$ is Lipschitz continuous with respect to $(t,y)$ on $J_T$.
\end{theo}

Our last main result is the following:

\begin{theo}\label{convergence2} {\bf (Convergence of the numerical solution under uniqueness assumption)}\\
Assume (A0)-(A1'). Let $T>0$ and $\varepsilon = (\Delta t, \Delta x)$ such that the CFL condition (\ref{CFL'}) is satisfied.  If $u :=\left( u^{\alpha} \right)_{\alpha}$ is the unique solution of (\ref{HJ})-(\ref{CI}) in the sense of Definition~\ref{defviscosity}, then the numerical solution $(U^{\alpha,n}_i)$ of (\ref{NS})-(\ref{CINS}) converges locally uniformly to $u$ when $\varepsilon$ goes to zero, on any compact set $\mathcal{K} \subset [0,T) \times J$, i.e.
\begin{equation}\label{eq::limit2}
\limsup_{\varepsilon\to 0} \ \sup_{(n \Delta t,i \Delta x) \in \mathcal{K}} \ | u^{\alpha} (n \Delta t,i \Delta x) - U^{\alpha, n}_i | = 0,
\end{equation}
where the index $\alpha$ in (\ref{eq::limit2}) is chosen such that $(n \Delta t,i \Delta x) \in \mathcal{K} \cap [0,T) \times J_{\alpha}$.
\end{theo}

Using our scheme (\ref{NS}), we will present in Section \ref{simulation} illustrations by numerical simulations with application to traffic.

\subsection{Brief review of the literature}

{\bf Hamilton-Jacobi formulation.} \ 
We mainly refer here to the comments provided in \cite{IMZ} and references there in. 
There is a huge literature dealing with HJ equations and mainly with equations with discontinuous Hamiltonians. 
However, concerning the study of HJ equation on a network, there exist a few works: 
the reader is referred to \cite{ACCT11,ACCT12} for a general definition of viscosity solutions on a network,
and \cite{CMS} for Eikonal equations. 
Notice that in those works, the Lagrangians depend on the position $x$ and are continuous with respect to this variable.
Conversely, in \cite{IMZ} the Lagrangians do not depend on the position but they are allowed to be discontinuous at the junction point.
Even for discontinuous Lagrangians, the uniqueness of the viscosity solution 
has been established in \cite{IMZ}.\\

{\bf Numerical schemes for Hamilton-Jacobi equations.} \ 
Up to our knowledge, there are no numerical schemes for HJ equations on junctions (except the very recent work \cite{GHZ}, see our Section \ref{sect::application} for more details),
while there are a lot of schemes for HJ equations for problems without junctions.
The majority of numerical schemes which were proposed to solve HJ equations are based on finite difference methods; 
see for instance \cite{CL} for upwind and centered discretizations, and \cite{FF,OSh} for ENO or WENO schemes. For finite elements methods, the reader could also refer to \cite{HS} and \cite{ZS}. Explicit classical monotone schemes have convergence properties 
but they require to satisfy a CFL condition and they exhibit a viscous behaviour. We can also cite Semi-Lagrangian schemes \cite{CD,Falcone,FF}.
Anti-diffusive methods coming from numerical schemes adapted for conservation laws were thus introduced \cite{BZ,XS}. 
Some other interesting numerical advances are done along the line of discontinuous Galerkin methods \cite{CS,BCS}.
Notice that more generally, an important effort deals with Hamilton-Jacobi-Bellman equations and Optimal Control viewpoint. It is out of the scope here.

\subsection{Organization of the paper}

In Section \ref{gradient}, we point out our first main property, namely Theorem~\ref{gradientestimate} about the time and space gradient estimates.
Then in Section \ref{s::convergence}, we first recall the notion of viscosity solutions for HJ equations. We then prove the second main property of our numerical scheme, namely Theorem~\ref{convergence1} and Theorem~\ref{convergence2} about the convergence of the numerical solution toward a solution of HJ equations when the mesh grid goes to zero.
In Section \ref{sect::application}, we propose the interpretation of our numerical results to traffic flows problems on a junction. In particular, the numerical scheme for HJ equations (\ref{NS}) is derived and the junction condition is interpreted. Indeed, we recover the well-known junction condition of Lebacque (see \cite{Lebacque96}) or equivalently those for the Riemann solver at the junction as in the book of Garavello and Piccoli \cite{GP}. 
Finally, in Section \ref{simulation} we illustrate the numerical behaviour of our scheme for a junction with two incoming and two outgoing branches.

\section{Gradient estimates for the scheme}
\label{gradient}

This section is devoted to the proofs of the first main result namely the time and space gradient estimates.

\subsection{Proof of Proposition~\ref{pro1}}

We begin by proving the monotonicity of the numerical scheme.

\begin{Proof}{\bf of Proposition \ref{pro1}:} We consider the numerical scheme given by (\ref{NS}) that we rewrite as follows for $n\ge 0$:
\begin{equation}\label{defSN1}
\begin{cases}
U^{\alpha ,n+1}_{i} =  S_{\alpha} \left[ U^{\alpha,n}_{i-1}, U^{\alpha,n}_i , U^{\alpha,n}_{i+1} \right] & \quad \text{for} \quad  i \ge 1, \quad \alpha=1,...,N ,\\
\\
U^{n+1}_{0} = S_0 \left[ U^{n}_0, (U^{\beta,n}_{1})_{\beta=1,...,N}  \right] & \quad \text{for} \quad i=0 ,
\end{cases}
\end{equation}
where
\begin{equation}\label{defSN2}
\begin{cases}
S_{\alpha} \left[ U^{\alpha,n}_{i-1}, U^{\alpha,n}_i , U^{\alpha,n}_{i+1} \right] &:= U^{\alpha ,n}_{i}- \Delta t \max \left\{ H^{+}_{\alpha} \left(\dfrac{U^{\alpha,n}_i - U^{\alpha,n}_{i-1}}{\Delta x} \right), H^{-}_{\alpha} \left( \dfrac{U^{\alpha,n}_{i+1} - U^{\alpha,n}_i}{\Delta x} \right) \right\},\\
\\
S_0 \left[ U^{n}_0, (U^{\beta,n}_{1})_{\beta=1,...,N}  \right]  &:= U^{n}_{0} - \Delta t \displaystyle{\max_{\beta = 1,\ldots ,N}} \ H^{-}_{\beta} \left( \dfrac{U^{\beta,n}_{1} - U^{n}_0}{\Delta x} \right) .
\end{cases}
\end{equation}

Checking the monotonicity of the scheme means checking that $S_{\alpha}$ and $S_0$ are non-decreasing in all their variables.\\
{\bf Case 1: $i\ge 1$}\\
This case is very classical. It is straightforward to check that $S_{\alpha}$ for any $\alpha=1,...,N$ is non-decreasing in 
$U^{\alpha,n}_{i-1}$ and $U^{\alpha,n}_{i+1}$. We compute
$$\frac{\partial S_{\alpha}}{\partial U^{\alpha,n}_i} = 
\left\{\begin{array}{ll}
\displaystyle 1 - \frac{\Delta t}{\Delta x} (H_\alpha^+)'(p_{i,-}^{\alpha,n}) & 
\quad \mbox{if}\quad \max \left\{ H_\alpha^+(p_{i,-}^{\alpha,n}), H_\alpha^-(p_{i,+}^{\alpha,n} \right\} = H_\alpha^+(p_{i,-}^{\alpha,n}),\\
\\
\displaystyle 1 - \frac{\Delta t}{\Delta x} (H_\alpha^-)'(p_{i,+}^{\alpha,n}) & 
\quad \mbox{if}\quad \max \left\{ H_\alpha^+(p_{i,-}^{\alpha,n}), H_\alpha^-(p_{i,+}^{\alpha,n} \right\} = H_\alpha^-(p_{i,+}^{\alpha,n})
\end{array}\right.$$
which is non-negative if the CFL condition  (\ref{CFL}) is satisfied.\\
{\bf Case 2: $i=0$}\\
Similarly, it is straightforward to check that $S_0$ is non-decreasing in each $U^{\beta,n}_{1}$ for $\beta=1,...,N$. We compute
$$\frac{\partial S_0}{\partial U^{n}_0} = 1 - \frac{\Delta t}{\Delta x} (H_\alpha^-)'(p_{0,+}^{\alpha,n}) \quad \mbox{if} \quad 
H_\alpha^-(p_{0,+}^{\alpha,n}) > H_\beta^-(p_{0,+}^{\beta,n}) \quad \mbox{for all}\quad  \beta \in \left\{1,...,N\right\}\backslash \left\{\alpha\right\}$$
which is also non-negative due to the CFL condition  (\ref{CFL}).\\
From cases 1 and 2, we deduce that the scheme is monotone.
\end{Proof}

\subsection{Proof of Theorem~\ref{gradientestimate}}

In this subsection, we prove the first main result Theorem~\ref{gradientestimate} about time and space gradient estimates.

Let us first define for any $n\geq 0$
\begin{equation}\label{defsupinfW}
m^n := \displaystyle{\inf_{\alpha,i}} \ W^{\alpha,n}_i \quad \mbox{and} \quad M^n := \displaystyle{\sup_{\alpha,i}} \ W^{\alpha,n}_i ,
\end{equation}
where $W^{\alpha,n}_i$ represents the time gradient defined in (\ref{defW}).

We also define
\begin{equation}\label{defsegmentI}
I^{\alpha,n}_{i,\sigma} := \begin{cases}
\left[p_{i,\sigma}^{\alpha,n},p_{i,\sigma}^{\alpha,n+1} \right] \quad \mbox{if} \quad p_{i,\sigma}^{\alpha,n} \leq p_{i,\sigma}^{\alpha,n+1},\\
\\
\left[p_{i,\sigma}^{\alpha,n+1},p_{i,\sigma}^{\alpha,n} \right] \quad \mbox{if} \quad p_{i,\sigma}^{\alpha,n} \geq p_{i,\sigma}^{\alpha,n+1}.
\end{cases} \quad \mbox{for} \quad \sigma \in \{ +,- \},
\end{equation}
with $p^{\alpha,n}_{i,\sigma}$ defined in (\ref{defp}) and we set
\begin{equation}\label{defD}
D^{\alpha,n}_{i,+}:=\displaystyle{\sup_{p_{\alpha}\in I^{\alpha,n}_{i,+}}} |H_{\alpha}'(p_{\alpha})|.
\end{equation}

In order to establish Theorem~\ref{gradientestimate}, we need the two following results namely Proposition~\ref{pro:0} and Lemma~\ref{lem:0}:

\begin{pro}\label{pro:0}{\bf (Time derivative estimate)}\\
Assume (A1).
Let $n \geq 0$ fixed and $\Delta x$, $\Delta t >0$. Let us consider $\left( U^{\alpha,n}_i \right)_{\alpha,i}$ satisfying for some constant $C^n>0$:
\begin{equation}\label{grad0}
|p^{\alpha,n}_{i,+}| \leq C^n, \quad \mbox{for} \quad i \geq 0, \quad \alpha=1,...N.
\end{equation}
We also consider $\left( U^{\alpha,n+1}_i \right)_{\alpha,i}$ and $\left( U^{\alpha,n+2}_i \right)_{\alpha,i}$ computed using the scheme (\ref{NS}).

If we have
\begin{equation}\label{CFLenD}
D^{\alpha,n}_{i,+} \leq \dfrac{\Delta x}{\Delta t} \quad \mbox{for any} \quad i \geq 0 \quad \mbox{and} \quad \alpha=1,...,N,
\end{equation}

Then it comes that
$$m^n \leq m^{n+1} \leq M^{n+1} \leq M^n.$$
\end{pro}

\begin{Proof}\\
{ \bf Step 0: Preliminaries. } \\
We introduce for any $n \geq 0$, $\alpha=1,...,N$ and for any $i \geq 1$, $\sigma \in \{+,-\}$ or for $i=0$ and $\sigma=+$:
\begin{equation}\label{defC}
C^{\alpha ,n}_{i,\sigma} :=- \sigma \displaystyle{\int^1_0} d\tau ( H^{-\sigma}_{\alpha})' (p^{\alpha ,n+1}_{i,\sigma} + \tau (p^{\alpha ,n}_{i,\sigma} - p^{\alpha ,n+1}_{i,\sigma})) \geq 0 .
\end{equation}

Notice that $C^{\alpha ,n}_{i,\sigma}$ is defined as the integral of $(H^{-\sigma}_{\alpha})'$ over a convex combination of $p$ with $p \in I^{\alpha,n}_{i,\sigma}$.
Hence for any $n \geq 0$, $\alpha=1,...,N$ and for any $i \geq 1$, $\sigma \in \{+,-\}$ or for $i=0$ and $\sigma=+$, we can check that
\begin{equation}\label{CFLforproof}
C^{\alpha,n}_{i,\sigma} \leq  \sup_{\substack{\beta=1,...,N \\ j \geq 0}} \ D^{\beta,n}_{j,+}.
\end{equation}

We also underline that for any $n \geq 0$, $\alpha=1,...,N$ and for any $i \geq 1$, $\sigma \in \{+,-\}$ or for $i=0$ and $\sigma=+$, we have the following relationship:
\begin{equation}\label{relationpW}
\dfrac{p^{\alpha,n}_{i,\sigma} - p^{\alpha,n+1}_{i,\sigma}}{\Delta t} = - \sigma \dfrac{W^{\alpha,n}_{i+\sigma} -W^{\alpha,n}_{i}}{\Delta x}.
\end{equation}

Let $n \geq 0$ be fixed and consider $\left( U^{\alpha,n}_i \right)_{\alpha,i}$ with $\Delta x$, $\Delta t >0$ given.
We compute $\left( U^{\alpha,n+1}_i \right)_{\alpha,i}$ and $\left( U^{\alpha,n+2}_i \right)_{\alpha,i}$ using the scheme (\ref{NS}).

{ \bf Step 1: Estimate on $m^n$ }\\
We want to show that $W^{\alpha,n+1}_i \geq m^n$ for any $i \geq 0$ and $\alpha=1,...,N$. It is then sufficient to take the infimum over $i \geq 0$ and $\alpha=1,...,N$ to conclude that
$$m^{n+1} \geq m^n.$$

Let $i \geq 0$ be fixed and we distinguish two cases:

{ \bf Case 1: Proof of $W^{\alpha,n+1}_{i} \geq m^{n} \quad \text{for all} \quad i \geq 1$}\\
Let a branch $\alpha$ fixed. We assume that
\begin{equation}\label{cas:0}
\displaystyle{\max} \left\{ H^{+}_{\alpha} (p^{\alpha,n+1}_{i,-}), H^{-}_{\alpha} (p^{\alpha,n+1}_{i,+}) \right\} = H^{-\sigma}_{\alpha} (p^{\alpha,n+1}_{i,\sigma}) \quad \text{for one} \quad \sigma \in \{ +,- \}.
\end{equation}

We have
$$ \begin{aligned}
\dfrac{W^{\alpha,n+1}_{i} -W^{\alpha,n}_{i}}{\Delta t} &= \dfrac{1}{\Delta t} \left( \displaystyle{\max} \left\{ H^{+}_{\alpha} (p^{\alpha,n}_{i,-}), H^{-}_{\alpha} (p^{\alpha,n}_{i,+}) \right\} - \displaystyle{\max} \left\{ H^{+}_{\alpha} (p^{\alpha,n+1}_{i,-}), H^{-}_{\alpha} (p^{\alpha,n+1}_{i,+}) \right\} \right) \\
&\geq \dfrac{1}{\Delta t} \left( H^{-\sigma}_{\alpha} (p^{\alpha,n}_{i,\sigma}) - H^{-\sigma}_{\alpha} (p^{\alpha,n+1}_{i,\sigma}) \right) \\
&= \dfrac{1}{\Delta t} \displaystyle{\int^1_0} d\tau ( H^{-\sigma}_{\alpha})' (p^{\alpha,n+1}_{i,\sigma} + \tau p) p \qquad \mbox{with} \quad p = p^{\alpha,n}_{i,\sigma} - p^{\alpha,n+1}_{i,\sigma} \\
&= C^{\alpha,n}_{i,\sigma} \displaystyle{\left(\dfrac{W^{\alpha,n}_{i+\sigma} -W^{\alpha,n}_{i}}{\Delta x}\right)}
\end{aligned} $$
where we use (\ref{relationpW}) and (\ref{defC}) in the last line.

Using (\ref{CFLforproof}) and (\ref{CFLenD}), we thus get
$$ \begin{aligned}
W^{\alpha,n+1}_{i} &\geq  \left(1 - C^{\alpha,n}_{i,\sigma} \dfrac{\Delta t}{\Delta x} \right) W^{\alpha,n}_{i} + C^{\alpha,n}_{i,\sigma} \dfrac{\Delta t}{\Delta x} W^{\alpha,n}_{i+\sigma} \\
&\geq  \displaystyle{\min(W^{\alpha,n}_{i} , W^{\alpha,n}_{i+\sigma})} \\
&\geq m^{n}.
\end{aligned} $$

{ \bf Case 2: Proof of $W^{n+1}_{i} \geq m^{n} \ \text{for} \ i=0$ }\\
We recall that in this case, we have $U^{\beta, n}_0 =: U^{n}_0$ for any $\beta=1,...,N$. Let us denote $W^{\beta,n}_{0}=:W^{n}_0=\dfrac{U^{n+1}_0-U^{n}_0}{\Delta t}$ for any $\beta=1,...,N$.
Then we define $\alpha_0$ such that
$$ H^{-}_{\alpha_0} (p^{\alpha_0,n+1}_{0,+}) = \displaystyle{\max_{\alpha=1,...,N} H^{-}_{\alpha} (p^{\alpha,n+1}_{0,+})}. $$

We argue like in Case 1 above and we get
$$ \dfrac{ W^{n+1}_0 - W^{n}_0 }{\Delta t} \geq C^{\alpha_0 ,n}_{0,+} \displaystyle{\left(\dfrac{W^{\alpha_0 ,n}_{1} -W^{n}_{0}}{\Delta x}\right)} .$$
Then using (\ref{CFLforproof}) and (\ref{CFLenD}) we conclude that:
$$W^{n+1}_0 \geq m^n.$$

{ \bf Step 2: : Estimate on $M^n$}\\
We recall that $n \geq 0$ is fixed. The proof for $M^n$ is directly adapted from Part 1.
We want to show that $W^{\alpha,n+1}_i \leq M^n$ for any $i \geq 0$ and $\alpha=1,...,N$. We distinguish the same two cases:
\begin{itemize}
\item[$\bullet$] If $i \geq 1$, instead of (\ref{cas:0}) we simply choose $\sigma$ such that
$$\displaystyle{\max} \left\{ H^{+}_{\alpha} (p^{\alpha,n}_{i,-}), H^{-}_{\alpha} (p^{\alpha,n}_{i,+}) \right\} = H^{-\sigma}_{\alpha} (p^{\alpha,n}_{i,\sigma}) \quad \text{for one} \quad \sigma \in \{ +,- \}.$$
\item[$\bullet$] If $i=0$, we define $\alpha_0$ such that
$$H^{-}_{\alpha_0} (p^{\alpha_0,n}_{0,+}) = \displaystyle{\max_{\alpha=1,...,N} H^{-}_{\alpha} (p^{\alpha,n}_{0,+})}. $$
\end{itemize}
Then taking the supremum, we can easily prove that
$$ M^{n+1} \leq M^{n}, \quad \mbox{for any} \quad n \geq 0.$$

By definition of $m^n$ and $M^n$ for a given $n \geq 0$, we recover the result
$${m}^{n} \leq  {m}^{n+1} \leq  {M}^{n+1}  \leq  {M}^{n}.$$
\end{Proof}

The second important result needed for the proof of Theorem~\ref{gradientestimate} is the following one:

\begin{lem}\label{lem:0}{\bf (Gradient estimate)}\\
Assume (A1).
Let $n \geq 0$ fixed and $\Delta x$, $\Delta t >0$. We consider that $\left( U^{\alpha,n}_i \right)_{\alpha,i}$ is given and we compute $\left( U^{\alpha,n+1}_i \right)_{\alpha,i}$ using the scheme (\ref{NS}).

If there exists a constant $K \in \R$ such that for any $i \geq 0$ and $\alpha=1,...,N$, we have
$$K \leq W^{\alpha,n}_{i}:=\dfrac{U^{\alpha,n+1}_{i}-U^{\alpha,n}_{i}}{\Delta t}$$
Then it follows for any $i \geq 0$ and $\alpha=1,...,N$
$$(H^{-}_{\alpha})^{-1}(-K) \leq p^{\alpha,n}_{i,+} \leq (H^{+}_{\alpha})^{-1}(-K)$$
with $p^{\alpha,n}_{i,+}$ defined in (\ref{defp}) and $(H^{-}_{\alpha})^{-1}$, $(H^{+}_{\alpha})^{-1}$ defined in (\ref{rem:1}).
\end{lem}

\begin{Proof}\\
Let $n \geq 0$ be fixed and consider $\left( U^{\alpha,n}_i \right)_{\alpha,i}$ with $\Delta x$, $\Delta t >0$ given.
We compute $\left( U^{\alpha,n+1}_i \right)_{\alpha,i}$ using the scheme (\ref{NS}).

Let us consider any $i \geq 0$ and $\alpha =1,...,N$. We distinguish two cases according to the value of $i$.

{ \bf Case 1: $i \geq 1$}\\
Assume that we have
$$ K \leq  W^{\alpha,n}_i=-\displaystyle{\max_{\sigma \in \{ +,- \}}} H_{\alpha}^{-\sigma} (p^{\alpha,n}_{i,\sigma}).$$
It is then obvious that we get
$$ -K \geq  H_{\alpha}^{-\sigma} (p^{\alpha,n}_{i,\sigma}), \quad \mbox{for any} \quad \sigma \in \{ +,- \}.$$

According to (A1) on the monotonicity of the Hamiltonians $H_{\alpha}$, we obtain
\begin{equation}\label{result1}
\begin{cases}
(H^{+}_{\alpha})^{-1} (-K) \geq p_{i,-}^{\alpha,n}=p_{i-1,+}^{\alpha,n} \\
\\
(H^{-}_{\alpha})^{-1} (-K) \leq  p_{i,+}^{\alpha,n}
\end{cases} \quad \mbox{for any} \quad i \geq 1, \quad n \geq 0 \quad \mbox{and} \quad \alpha =1, \ldots , N.
\end{equation}

{ \bf Case 2: $i=0$}\\
The proof is similar to Case 1 because on the one hand we have
$$ K \leq  W^{\alpha,n}_0 =: W^{n}_0 = -\displaystyle{\max_{\beta=1,...,N}} H_{\beta}^{-} (p^{\beta,n}_{0,+}) $$
which obviously leads to
$$ (H_{\alpha}^{-})^{-1}(-K) \leq p^{\alpha,n}_{0,+},$$
where we use the monotonicity of $H_{\alpha}^{-}$ from assumption (A1).
On the other hand, from (\ref{result1}) we get 
$$ (H^{+}_{\alpha})^{-1}(-K) \geq p_{1,-}^{\alpha,n}=p_{0,+}^{\alpha,n}.$$

We conclude
$$ (H^{-}_{\alpha})^{-1} (-K) \leq p^{\alpha,n}_{i,+} \leq (H^{+}_{\alpha})^{-1} (-K), \quad \mbox{for any} \quad i, \quad n \geq 0 \quad \mbox{and} \quad \alpha =1,...,N $$
which ends the proof.
\end{Proof}

\begin{Proof}{\bf of Theorem~\ref{gradientestimate}:}
The idea of the proof is to introduce new continuous Hamiltonians $\tilde{H}_{\alpha}$ that satisfy the following properties:
\begin{enumerate}[(i)]
\item the new Hamiltonians $\tilde{H}_{\alpha}$ are equal to the old ones $H^{\alpha}$ on the segment $[\underline{p}_{\alpha},\overline{p}_{\alpha}]$,
\item the derivative of the new Hamiltonians $|\tilde{H}_{\alpha}'|$ taken at any point is less or equal to $\displaystyle \sup_{p \in [\underline{p}_{\alpha},\overline{p}_{\alpha}]} |H_{\alpha}'(p)|$.
\end{enumerate}
This modification if the Hamiltonians is done in order to show that the gradient stays in the interval $[\underline{p}_{\alpha},\overline{p}_{\alpha}]$.

{ \bf Step 1: Modification of the Hamiltonians}\\
Let the new Hamiltonians $\tilde{H}_{\alpha}$ for all $\alpha =1,...,N$ be defined as
\begin{equation}\label{defHtilde}
\tilde{H}_{\alpha}(p)=
\begin{cases}
g_{\alpha}^l(p) \quad &\text{for} \quad p \leq \underline{p}_{\alpha} \\
H_{\alpha}(p) \quad &\text{for} \quad p \in [\underline{p}_{\alpha},\overline{p}_{\alpha}] \\
g_{\alpha}^r(p) \quad &\text{for} \quad p \geq \overline{p}_{\alpha}
\end{cases}
\end{equation}
with $\underline{p}_{\alpha}$, $\overline{p}_{\alpha}$ defined in (\ref{defcontgrad}) and $g_{\alpha}^l$, $g_{\alpha}^r$ two functions such that
$$\begin{cases}
g_{\alpha}^l \in C^1((-\infty,\underline{p}_{\alpha}]), \\
g_{\alpha}^l(\underline{p}_{\alpha})=-m_0, \\
(g_{\alpha}^l)'(\underline{p}_{\alpha})=H'_{\alpha}(\underline{p}_{\alpha}), \\
(g_{\alpha}^l)'<0 \quad &\mbox{on} \quad (-\infty,\underline{p}_{\alpha}), \\
|(g_{\alpha}^l)'(p)|< |H'_{\alpha}(\underline{p}_{\alpha})| \quad &\mbox{for} \quad p<\underline{p}_{\alpha}, \\
g_{\alpha}^l \to +\infty \quad &\mbox{for} \quad p \to -\infty,
\end{cases} 
\quad \mbox{and} \quad 
\begin{cases}
g_{\alpha}^r \in C^1([\overline{p}_{\alpha},+\infty)), \\
g_{\alpha}^r(\overline{p}_{\alpha})=-m_0, \\
(g_{\alpha}^r)'(\overline{p}_{\alpha})=H'_{\alpha}(\overline{p}_{\alpha}), \\
(g_{\alpha}^r)'>0 \quad &\mbox{on} \quad (\overline{p}_{\alpha},+\infty), \\
|(g_{\alpha}^r)'(p)|< |H'_{\alpha}(\overline{p}_{\alpha})| \quad &\mbox{for} \quad p>\overline{p}_{\alpha}, \\
g_{\alpha}^r \to +\infty \quad &\mbox{for} \quad p \to +\infty.
\end{cases}$$ 

We can easily check that
\begin{equation}\label{cond1}
0< \tilde{H}'_{\alpha} <\displaystyle{\sup_{p_{\alpha}\in [\underline{p}_{\alpha},\overline{p}_{\alpha}]}} |H_{\alpha}'(p_{\alpha})|, \quad \mbox{on} \quad \R \setminus [\underline{p}_{\alpha},\overline{p}_{\alpha}],
\end{equation}
and
\begin{equation}\label{cond2}
\tilde{H}_{\alpha}>-m_0 \quad \mbox{on} \quad \R \setminus [\underline{p}_{\alpha},\overline{p}_{\alpha}].
\end{equation}

We can also check that $\tilde{H}_{\alpha}$ satisfies (A1). Then Proposition \ref{pro:0} and Lemma \ref{lem:0} hold true for the new Hamiltonians $\tilde{H}_{\alpha}$ (especially we can adapt (\ref{rem:1}) to the $\tilde{H}_{\alpha}$ for defining a sort of inverse).

Let $\tilde{H}^{+}_{\alpha}$ (resp. $\tilde{H}^{-}_{\alpha}$) denotes the non-decreasing (resp. non-increasing) part of $\tilde{H}_{\alpha}$.

We consider the new numerical scheme for any $n \geq 0$ that reads as:
\begin{equation}\label{NS'}
\begin{cases}
\dfrac{\tilde{U}^{\alpha ,n+1}_{i} -\tilde{U}^{\alpha ,n}_{i}}{\Delta t} + \max \left\{ \tilde{H}^{+}_{\alpha} (\tilde{p}^{\alpha ,n}_{i,-}), \tilde{H}^{-}_{\alpha} (\tilde{p}^{\alpha ,n}_{i,+}) \right\}=0, \quad &\mbox{for}\quad i \ge 1, \quad \alpha = 1,\ldots,N ,\\
\\
\left.\begin{array}{l}
\tilde{U}^{\beta,n}_{0} =: \tilde{U}^{n}_{0}, \quad \text{for all} \quad \beta = 1,\ldots,N \\
\\
\dfrac{\tilde{U}^{n+1}_{0} -\tilde{U}^{n}_{0}}{\Delta t} + \displaystyle{\max_{\beta = 1,\ldots ,N}} \ \tilde{H}^{-}_{\beta} (p^{\beta ,n}_{0,+}) =0
\end{array}\right| \quad &\mbox{for} \quad i=0,
\end{cases}
\end{equation}
subject to the initial condition
\begin{equation}\label{CINS'}
\tilde{U}^{\alpha,0}_{i} = U^{\alpha,0}_{i} = u^{\alpha}_{0} (i \Delta x), \quad i \ge 0, \quad \alpha = 1,\ldots,N .
\end{equation}

The discrete time and space gradients are defined such as:
\begin{equation}\label{defW'}
\tilde{p}^{\alpha ,n}_{i,+}:=\dfrac{\tilde{U}^{\alpha ,n}_{i+1} -\tilde{U}^{\alpha ,n}_{i}}{\Delta x} \quad \text{and} \quad \tilde{W}^{\alpha,n}_i := \dfrac{\tilde{U}^{\alpha ,n+1}_{i} -\tilde{U}^{\alpha ,n}_{i}}{\Delta t} .
\end{equation}

Let us consider 
\begin{equation}\label{defsupinfW'}
\tilde{m}^n := \displaystyle{\inf_{i,\alpha}} \ \tilde{W}^{\alpha,n}_{i} \quad \mbox{and} \quad \tilde{M}^n := \displaystyle{\sup_{i,\alpha}} \ \tilde{W}^{\alpha,n}_{i}
\end{equation}
where $\tilde{W}^{\alpha,n}_{i}$ is defined in (\ref{defW'}).
We also set
\begin{equation}\label{defD2}
\tilde{D}^{\alpha,n}_{i,+}:=\displaystyle{\sup_{p_{\alpha}\in \tilde{I}^{\alpha,n}_{i,+}}} |\tilde{H}_{\alpha}'(p_{\alpha})|,
\end{equation}
where $\tilde{I}^{\alpha,n}_{i,+}$ is the analogue of $I^{\alpha,n}_{i,+}$ defined in (\ref{defsegmentI}) with $\tilde{p}^{\alpha ,n}_{i,+}$ and $\tilde{p}^{\alpha ,n+1}_{i,+}$ given in (\ref{defW'}).

According to (\ref{cond1}), the supremum of $|\tilde{H}'_{\alpha}|$ is reached on $[\underline{p}_{\alpha},\overline{p}_{\alpha}]$. As $\tilde{H}_{\alpha} \equiv H_{\alpha}$ on $[\underline{p}_{\alpha},\overline{p}_{\alpha}]$, the CFL condition~(\ref{CFL'}) gives that for any $i \geq 0$, $n \geq 0$ and $\alpha=1,...,N$:
\begin{equation}\label{CFLproof}
\tilde{D}^{\alpha,n}_{i,+} \leq \sup_{p_{\alpha}\in [\underline{p}_{\alpha},\overline{p}_{\alpha}]} |H_{\alpha}'(p_{\alpha})| \leq \dfrac{\Delta x}{\Delta t}.
\end{equation}

{ \bf Step 2: First gradient bounds}\\
Let $n \geq 0$ be fixed. By definition~(\ref{defsupinfW'}) and if $\tilde{m}^n$ is finite, we have
$$\tilde{m}^n \leq \tilde{W}^{\alpha,n}_i, \quad \mbox{for any} \quad i \geq 0, \quad \alpha=1,...,N.$$
Using Lemma~\ref{lem:0}, it follows that
\begin{equation}
(\tilde{H}^-_{\alpha})^{-1}(-\tilde{m}^n) \leq \tilde{p}^{\alpha,n}_{i,+} \leq  (\tilde{H}^+_{\alpha})^{-1}(-\tilde{m}^n), \quad \mbox{for any} \quad i \geq 0 \quad \mbox{and} \quad \alpha=1,...,N.
\end{equation}
We define 
$$C^n=\max \left\{ \left|(\tilde{H}^-_{\alpha})^{-1}(-\tilde{m}^n) \right|, \left|(\tilde{H}^+_{\alpha})^{-1}(-\tilde{m}^n) \right| \right\} >0,$$
and we recover that
$$|\tilde{p}^{\alpha,n}_{i,+}| \leq C^n, \quad \mbox{for any} \quad i \geq 0, \quad \alpha=1,...,N.$$

{ \bf Step 3: Time derivative and gradient estimates}\\
For any $n \geq 0$, (\ref{CFLproof}) holds true. Moreover, if $\tilde{m}^n$ is finite, then there exists $C^n >0$ such that
$$|\tilde{p}^{\alpha,n}_{i,+}| \leq C^n, \quad \mbox{for any} \quad i \geq 0, \quad \alpha=1,...,N.$$

Then using Proposition~\ref{pro:0} we get
\begin{equation}\label{timest'}
\tilde{m}^n \leq \tilde{m}^{n+1} \leq \tilde{M}^{n+1} \leq \tilde{M}^{n} \quad \mbox{for any} \quad n \geq 0.
\end{equation}
In particular, $\tilde{m}^{n+1}$ is also finite.

Using the assumption that $m^0$ is finite and according to (\ref{defcontgrad}), Lemma \ref{lem:0} and the scheme (\ref{NS}), we can check that
\begin{equation}\label{encadrement}
\underline{p}_{\alpha} \leq p^{\alpha,0}_{i,+} \leq \overline{p}_{\alpha} \quad \mbox{for any} \quad i \geq 0 \quad \mbox{and} \quad \alpha=1,...,N.
\end{equation}

From (\ref{CINS'}), we have $p^{\alpha,0}_{i,+}=\tilde{p}^{\alpha,0}_{i,+}$. Therefore, from (\ref{defHtilde}), (\ref{NS}) and (\ref{encadrement}), we deduce that $\tilde{W}^{\alpha,0}_i=W^{\alpha,0}_i$ and we obtain that
$$\tilde{m}^0=m^0.$$

According to (\ref{timest'}), we deduce that $m^0 \leq \tilde{W}^{\alpha,n}_i$ for any $i\geq 0$, $n \geq 0$ and $\alpha=1,...,N$. 

Then using Lemma~\ref{lem:0} and (\ref{cond2}), we conclude that for all $i\geq 0$, $n \geq 0$ and $\alpha=1,...,N$
\begin{equation}\label{gradest'}
\underline{p}_{\alpha} \leq \tilde{p}^{\alpha,n}_{i,+} \leq \overline{p}_{\alpha}.
\end{equation}

{ \bf Step 4: Conclusion}\\
If (\ref{gradest'}) holds true, then $\tilde{H}_{\alpha}(\tilde{p}^{\alpha,n}_{i,+})= H_{\alpha}(\tilde{p}^{\alpha,n}_{i,+})$ for all $i\geq 0$, $n \geq 0$ and $\alpha=1,...,N$. Thus the modified scheme (\ref{NS'}) is strictly equivalent to the original scheme (\ref{NS}) and $U^{\alpha,n}_{i}=\tilde{U}^{\alpha,n}_{i}$. We finally recover the results for all $i\geq 0$, $n \geq 0$ and $\alpha=1,...,N$:
\begin{enumerate}[(i)]
\item (Time derivative estimate)
$$m^n \leq m^{n+1} \leq M^{n+1} \leq M^{n},$$
\item (Gradient estimate)
$$\underline{p}_{\alpha} \leq p^{\alpha,n}_{i,+} \leq \overline{p}_{\alpha}.$$
\end{enumerate}
\end{Proof}

\begin{rem}{ \bf (Do the bounds (\ref{eq::g2}) always give informations on the gradient?)}\\
We assume that the Theorem \ref{gradientestimate} holds true.
\begin{enumerate}[(i)]
\item {\bf (Bounds on $m^n$)} From the scheme (\ref{NS}), we can rewrite
$$m^n=\inf_{\alpha,i} \ \min_{ \sigma \in \{ +,- \} } \left\{-H^{-\sigma}_{\alpha}(p^{\alpha,n}_{i,\sigma})\right\}.$$
It is then obvious that
$$-m^0 \geq \min_{p_{\alpha} \in \R} \ H_{\alpha}(p_{\alpha}) \quad \text{for} \quad \alpha=1,...,N,$$
which ensures that the bound from below in (\ref{eq::g2}) always gives an information on the gradient $(p^{\alpha,n}_{i,+})$.
\item {\bf (Bounds on $M^n$)} For the bounds from above in (\ref{eq::g2}), we get 
\begin{equation}\label{eq:1}
H_{\alpha} \left(p^{\alpha,n}_{i,+}\right) \geq -M^0 \quad \mbox{for all} \quad \alpha=1,...,N, \quad i \geq 0 \quad \mbox{and} \quad n \leq n_T.
\end{equation}
Note that for each $\alpha=1,...,N$, (\ref{eq:1}) gives an information on the $(p^{\alpha,n}_{i,+})$ only if
$$-M_0 > \min_{p_{\alpha} \in \R} H_{\alpha}(p_{\alpha}).$$
\end{enumerate}
\end{rem}

\begin{rem}\label{weakerassumptions}{\bf (Extension to weaker assumptions on  $H_{\alpha}$ than (A1))}\\
All the results of this paper can be extended if we consider weaker conditions than (A1) on the Hamiltonians $H_{\alpha}$. Indeed, we can assume that the $H_{\alpha}$ for any $\alpha=1,...,N$ are locally Lipschitz. 
This assumption is more adapted for our traffic application (see Section~\ref{sect::application}).

We now focus on what should be modified if we do so. 

{\bf How to modify CFL condition (\ref{CFL})?} \\
The main new idea is then to consider the closed convex hull for the discrete gradient defined by
$$I^{\alpha,n}:=\overline{Conv(p^{\alpha,n}_{i,+})}_{i \geq 0}.$$
Then the Lipschitz constant $L^{\alpha,n}$ of the considered $H_{\alpha}$ is a natural upper bound 
$$| H_{\alpha}(p+q)-H_{\alpha}(p) | \leq L^{\alpha,n} |q| \quad \mbox{with} \quad p, \ p+q \in I^{\alpha,n}.$$
Then the natural CFL condition which replaces (\ref{CFL}) is the following one:
\begin{equation}\label{CFLlipschitz}
\dfrac{\Delta x}{\Delta t} \geq \sup_{\substack{\alpha=1,...,N \\ n \leq n_T}} L^{\alpha,n}.
\end{equation}
With such a condition, we can easily prove the monotonicity of the numerical scheme.\\

{\bf How to modify CFL condition (\ref{CFL'})?} \\
Assume that CFL condition (\ref{CFL'}) is replaced by the following one
\begin{equation}\label{CFL'mod}
\dfrac{\Delta x}{\Delta t} \geq \displaystyle{\esssup_{\substack{\alpha=1,\dots ,N \\ p_{\alpha} \in [\underline{p}_{\alpha}, \overline{p}_{\alpha}]}}} |H'_{\alpha}(p_{\alpha}) |,
\end{equation}
where $\esssup$ denotes the essential supremum.

In the proof of Theorem~\ref{gradientestimate}, the time derivative estimate uses the integral of $H'_{\alpha}$ which is defined almost everywhere if $H_{\alpha}$ is at least Lipschitz. The remaining of the main results of Section \ref{mainresults} do not use a definition of $H'_{\alpha}$, except in the CFL condition. We just need to satisfy the new CFL condition~(\ref{CFL'mod}).
\end{rem}

\section{Convergence result for the scheme}
\label{s::convergence}

\subsection{Viscosity solutions}
\label{viscositysection}

We introduce the main definitions related to viscosity solutions for HJ equations that are used in the remaining.
For a more general introduction to viscosity solutions, the reader could refer to Barles \cite{B} and to Crandall, Ishii, Lions \cite{CIL}.

Let $T>0$. We set $u:=(u^{\alpha})_{\alpha=1,...,N} \in C^1_{*}(J_T)$ where $C^1_{*}(J_T)$ is defined in (\ref{spacefunct}) and we consider the additional condition
$$u^{\alpha}(t,0)=u^{\beta}(t,0)=:u(t,0) \quad \mbox{for any} \quad \alpha, \beta .$$

\begin{rem}
Following \cite{IMZ}, we recall that (\ref{HJ}) can be rigorously rewritten as
\begin{equation}\label{HJ2}
u_t + H(y,u_y)= 0, \quad \mbox{for} \quad (t,y) \in [0,T) \times J,
\end{equation}
with
$$H(y,p):=
\left\{ \begin{array}{lll}
H_{\alpha}(p), \quad &\mbox{for} \quad p \in \R, \quad &\mbox{if} \quad y \in J^*_{\alpha}, \\
\\
\displaystyle{\max_{\alpha=1,...,N}} H^-_{\alpha}(p_{\alpha}), \quad &\mbox{for} \quad p=(p_1,...,p_N) \in \R^{\N}, \quad &\mbox{if} \quad y=0,
\end{array} \right.$$
subject to the initial condition
\begin{equation}\label{CI2}
u(0,y)=u_0(y):=\left( u^{\alpha}_0(x) \right)_{\alpha=1,...,N}, \quad \mbox{for} \quad y=x e_{\alpha} \in J \quad \mbox{with} \quad x \in [0,+\infty).
\end{equation}
\end{rem}

\begin{defi}{\bf (Upper and lower semi-continuous envelopes)}\\
For any function $u : [0,T) \times J \to \R$, upper and lower semi-continuous envelopes are respectively defined as:
$$ u^* (t,y)=\limsup_{(t',y')\to (t,y)} \ u(t',y') \quad \text{and} \quad u_* (t,y)=\liminf_{(t',y')\to (t,y)} \ u(t',y') .$$

Moreover, we recall
$$ \left\{ \begin{array}{l}
u \quad \mbox{is upper semi-continuous if and only if} \quad  u=u^* ,\\
\\
u \quad  \mbox{is lower semi-continuous if and only if} \quad  u=u_* ,\\
\\
u \quad  \mbox{is continuous if and only if} \quad  u_*=u^* .
\end{array} \right. $$
\end{defi}

\begin{defi}\label{defviscosity}{\bf (Viscosity solutions)}\\
A function $u: [0,T) \times J \to \R$ is a viscosity subsolution (resp. supersolution) of (\ref{HJ}) on $J_T=(0,T) \times J$ if it is an upper semi-continuous (resp. lower semi-continuous) function, and if for any $P=(t,y)\in J_T$ and any test function $\varphi:=(\varphi^\alpha)_{\alpha} \in C^1_{*}(J_T)$ such that $u-\varphi \leq 0$ (resp. $u-\varphi \geq 0$) at the point $P$, we have
\begin{equation}\label{viscosityineq1A}
\varphi^{\alpha}_t (t,x) + H_{\alpha}(\varphi^{\alpha}_x(t,x)) \le 0 \quad \text{if} \quad y=x e_{\alpha} \in J^*_{\alpha},
\end{equation}

\begin{equation}\label{viscosityineq2A}
\Big( \text{resp.} \quad
\varphi^{\alpha}_t (t,x) + H_{\alpha}(\varphi^{\alpha}_x(t,x)) \ge 0 \quad \text{if} \quad y=x e_{\alpha} \in J^*_{\alpha}
\Big),
\end{equation}

\begin{equation}\label{viscosityineq1B}
\varphi_t (t,0) + \displaystyle{\max_{\alpha = 1,\ldots,N}} \ H^{-}_{\alpha}(\varphi^{\alpha}_x(t,0)) \le 0 \quad \text{if} \quad y=0 ,
\end{equation}

\begin{equation}\label{viscosityineq2B}
\Big( \text{resp.} \quad
\varphi_t (t,0) + \displaystyle{\max_{\alpha = 1,\ldots,N}} \ H^{-}_{\alpha}(\varphi^{\alpha}_x(t,0)) \ge 0 \quad \text{if} \quad y=0
\Big).
\end{equation}

A function $u^*$ (resp. $u_*$) is said to be a viscosity subsolution (resp. supersolution) of (\ref{HJ})-(\ref{CI}) on $[0,T) \times J$, if $u^*$ is a viscosity subsolution (resp. $u_*$ is a viscosity supersolution) of (\ref{HJ}) on $J_T$ and if moreover it satisfies:
$$ \begin{aligned}
\begin{cases}
u^*(0,y)\le u_0(y)  \quad &\text{for all} \quad  y \in J ,\\
\\
\Big( \mbox{resp.} \quad u_*(0,y) \ge u_0(y)  \quad &\text{for all} \quad y \in J \Big),
\end{cases}
\end{aligned} $$ 
when the initial data $u_0$ is assumed to be continuous.\\

A function $u: [0,T) \times J \to \R$ is said to be a viscosity solution of (\ref{HJ}) on $J_T$ (resp. of (\ref{HJ})-(\ref{CI}) on $[0,T) \times J$) if $u^*$ is a viscosity subsolution and $u_*$ is a viscosity supersolution of (\ref{HJ}) on $J_T$ (resp. of (\ref{HJ})-(\ref{CI}) on $[0,T) \times J$).
\end{defi}

Hereafter, we recall two properties of viscosity solutions on a junction that are extracted from~\cite{IMZ}:
\begin{pro}\label{pro:1}{\bf (Comparison principle)}\\
Assume (A0)-(A1') and let $T>0$. Assume that $\overline{u}$ and $\underline{u}$ are respectively a viscosity subsolution and a viscosity supersolution of (\ref{HJ})-(\ref{CI}) on $[0,T) \times J$ in the sense of Definition~\ref{defviscosity}.
We also assume that there exists a constant $C_T > 0$ such that for all $(t,y) \in [0,T) \times J$
$$\overline{u}(t,y) \leq C_T(1+|y|) \quad \left( \mbox{resp.} \quad \underline{u}(t,y) \geq -C_T(1+|y|) \right).$$
Then we have $\overline{u} \le \underline{u}$ on $J_T$.
\end{pro}

\begin{pro}\label{pro:4}{\bf (Equivalence with relaxed junction conditions)}\\
Assume (A1') and let $T>0$. A function $u :[0,T) \times J \to \R$ is a viscosity subsolution (resp. a viscosity supersolution) of (\ref{HJ}) on $J_T$ if and only if for any function $\varphi:=(\varphi^\alpha)_{\alpha} \in C^1_{*}(J_T)$ and for any $P=(t,y)\in J_T$ such that  $u-\varphi \leq 0$ (resp. $u-\varphi \geq 0$) at the point $P$, we have the following properties
\begin{itemize}
\item[$\bullet$] if $y=x e_{\alpha} \in J_{\alpha}^{*}$, then
$$ \varphi^{\alpha}_t(t,x)+H_{\alpha}(\varphi^{\alpha}_x(t,x)) \leq 0 \quad \text{(resp.} \geq 0 \text{)} $$
\item[$\bullet$] if $y=0$, then either there exists one index $\alpha \in \{1,...,N \}$ such that
$$ \varphi^{\alpha}_t(t,0)+H_{\alpha}(\varphi^{\alpha}_x(t,0)) \leq 0 \quad \text{(resp.} \geq 0 \text{)} $$
or (\ref{viscosityineq1B}) (resp. (\ref{viscosityineq2B})) holds true for $y=0$.
\end{itemize}
\end{pro}

We skip the proof of Proposition~\ref{pro:1} and Proposition~\ref{pro:4} which are directly available in \cite{IMZ}.

\subsection{Convergence of the numerical solution}
\label{sub::convergence}

We assume (A0), (A1') and we set $\varepsilon :=(\Delta t,\Delta x)$ satisfying the CFL condition (\ref{CFL'}). This section first deals with a technical result (see Lemma~\ref{lem:1}) that is very useful for the proof of Theorem~\ref{convergence2} that is the convergence of the numerical solution of (\ref{NS})-(\ref{CINS}) towards a solution of (\ref{HJ})-(\ref{CI}) when $\varepsilon$ goes to zero. According to Theorem~\ref{existence&uniqueness}, we know that the equation (\ref{HJ})-(\ref{CI}) admits a unique solution in the sense of Definition~\ref{defviscosity}.
For Theorem~\ref{convergence1}, we extend the convergence proof, assuming the weakest assumption (A1) instead of (A1').
We close this subsection with the proof of gradient estimates for the continuous solution (see Corollary~\ref{cor2}).\\

We denote by 
$$u^{\alpha}_{\varepsilon}(n \Delta t, i \Delta x):=U^{\alpha,n}_i$$
an approximation of $ u^{\alpha}(n \Delta t,i \Delta x)$ for any $\alpha=1,...,N$ and $i\geq 0$ , $n\geq 0$.
Assume that $u^{\alpha}_{\varepsilon}$ solves the numerical scheme (\ref{NS})-(\ref{CINS}).
We recall
$$u^{\alpha}_{\varepsilon}(n \Delta t,0)=: u_{\varepsilon}(n \Delta t,0), \quad \mbox{for any} \quad \alpha=1,...,N.$$

We also denote by $G^{\alpha}_{\varepsilon}$ a set of all grid points $(n \Delta t, i \Delta x)$ on $[0,T) \times J_{\alpha}$ for any branch $\alpha=1,...,N$, and we set 
\begin{equation}\label{grid}
G_{\varepsilon}=\displaystyle{\bigcup_{\alpha=1,...,N}} G^{\alpha}_{\varepsilon}
\end{equation}
the whole set of grid points on $[0,T) \times J$, with identification of the junction points $(n \Delta t,0)$ of each grid $G^{\alpha}_{\varepsilon}$.

We call $u_{\varepsilon}$ the function defined by its restrictions to the grid points of the branches
$$u_{\varepsilon}=u^{\alpha}_{\varepsilon} \quad \mbox{on} \quad G^{\alpha}_{\varepsilon}.$$

For any point $(t,y) \in [0,T) \times J$, we define the half relaxed limits
\begin{equation}\label{subsolution}
\overline{u}(t,y) = \limsup_{\substack{\varepsilon \to 0 \\  G_{\varepsilon}\ni (t',y') \to (t,y)}} \ u_{\varepsilon}(t',y'),
\end{equation}

\begin{equation}\label{supersolution}
\left( \mbox{resp.} \quad \underline{u}(t,y) = \liminf_{\substack{\varepsilon\to 0 \\ G_{\varepsilon}\ni (t',y') \to (t,y)}} \ u_{\varepsilon}(t',y') \right).
\end{equation}

Thus we have that $\overline{u}:=\left( \overline{u}^{\alpha} \right)_{\alpha}$ (resp. $\underline{u}:=\left( \underline{u}^{\alpha} \right)_{\alpha}$) is upper semi-continuous (resp. lower semi-continuous).\\


\begin{lem}\label{lem:1}{\bf ($\varepsilon$-uniform space and time gradient bounds)}\\
Assume (A0)-(A1). Let $T>0$ and $\varepsilon = (\Delta t, \Delta x)$ such that the CFL condition (\ref{CFL'}) is satisfied.
If $(U^{\alpha,n}_i)$ is the numerical solution of (\ref{NS})-(\ref{CINS}),
then for any $\alpha=1,...,N$, $i \geq 0$ and $n \geq 0$, we have
\begin{equation}\label{spacetimegrad}
\underline{p}^0_{\alpha} \leq \dfrac{U^{\alpha,n}_{i+1}-U^{\alpha,n}_{i}}{\Delta x} \leq \overline{p}^0_{\alpha} \quad \mbox{and} \quad m^0_0 \leq \dfrac{U^{\alpha,n+1}_{i}-U^{\alpha,n}_{i}}{\Delta t} \leq M^0_0,
\end{equation}
where the quantities (independent of $\varepsilon$) $m^0_0$, $M^0_0$, $\underline{p}^0_{\alpha}$ and $\overline{p}^0_{\alpha}$ are respectively defined in (\ref{defcontgrad1}) and (\ref{defcontgrad2}).
\end{lem}

\begin{Proof}{\bf of Lemma~\ref{lem:1}:} \\
Let $\varepsilon=(\Delta t, \Delta x)$ be fixed such that the CFL condition (\ref{CFL'}) is satisfied.\\
{ \bf Step 1: Proof of $m^0 \geq m^0_0$, $\underline{p}^0_{\alpha} \leq \underline{p}_{\alpha}$ and $\overline{p}_{\alpha} \leq \overline{p}^0_{\alpha}$}\\
We first show that
\begin{equation}\label{*10}
m^0 \geq m^0_0.
\end{equation} 
Indeed using (A1) and the fact that $H_{\alpha}(p) = \max \left\{H^-_{\alpha}(p),H^+_{\alpha}(p) \right\}$ for any $p \in \R$, we get
$$m^0 =\displaystyle{\inf_{\alpha,i}} \left\{ -H_{\alpha}(p^{\alpha,0}_i) \right\} \geq \displaystyle{\inf_{\substack{\alpha \\ p_{\alpha} \in [L^{\alpha,-}, L^{\alpha,+}]}}} \left\{ -H_{\alpha}(p_{\alpha}) \right\} =: m^0_0,$$
where we recall that $L^{\alpha,-}$ and $L^{\alpha,+}$ are the best Lipschitz constants defined in (\ref{LipConst}) that implies
\begin{equation}\label{*5}
L^{\alpha,-} \leq p^{\alpha,0}_{i,+} \leq L^{\alpha,+}, \quad \mbox{for any} \quad i \geq 0.
\end{equation}

From (\ref{defcontgrad2}) and the monotonicity of $H^{\pm}_{\alpha}$, we deduce
\begin{equation}\label{*6}
\underline{p}^0_{\alpha} \leq \underline{p}_{\alpha} \quad \mbox{and} \quad \overline{p}_{\alpha} \leq \overline{p}^0_{\alpha}, \quad \mbox{ for any} \quad \alpha=1,...,N.
\end{equation}

{\bf Step 2: Proof of $M^0 \leq M^0_0$}\\
Recall the definitions
$$M^0:=\sup_{\alpha,i} \ W^{\alpha,0}_i = \max \{A, B\}, \quad \mbox{with} \quad \begin{cases}
A:=\displaystyle{\min_{\alpha=1,...,N}} \left\{-H^{-}_{\alpha} (p^{\alpha,0}_{0,+}) \right\}, \\
B:=\displaystyle{\sup_{\alpha,i \geq 1}} \left( \displaystyle{\min_{\sigma \in \{ +,- \}}} \left\{-H^{-\sigma}_{\alpha} (p^{\alpha,0}_{i,\sigma}) \right\} \right).
\end{cases}$$
and
$$M^0_0:= \displaystyle{\max \left[ \max_{\alpha=1,...,N} \left(\min_{\sigma \in \{ +,- \}} \left\{ -H^{-\sigma}_{\alpha}(L^{\alpha,\sigma}) \right\} \right) , \min_{\alpha=1,...,N} \left\{ -H^{-}_{\alpha}(L^{\alpha,+}) \right\} \right]}.$$

Let us show that 
\begin{equation}\label{*11}
M^0 \leq M^0_0.
\end{equation}
We distinguish two cases according to the value of $M^0$:
\begin{itemize}
\item[$\bullet$] If $M^0=A$, then
$$M^0_0 \geq \displaystyle{\min_{\alpha=1,...,N}} (-H^-_{\alpha}(L^{\alpha,+})) \geq \displaystyle{\min_{\alpha=1,...,N}} (-H^-_{\alpha}(p^{\alpha,0}_{0,+})) =A= M^0,$$
where we use (\ref{*5}) and the monotonicity of $H^-_{\alpha}$.
\item[$\bullet$] If $M^0=B$, then
$$M^0_0 \geq \max_{\alpha=1,...,N} \ \left(\min_{\sigma \in \{ +,- \}} \ (-H^{-\sigma}_{\alpha}(L^{\alpha,\sigma})) \right) \geq \displaystyle{\sup_{\alpha,i \geq 1}} \left( \displaystyle{\min_{\sigma \in \{ +,- \}}} (-H^{-\sigma}_{\alpha} (p^{\alpha,0}_{i,\sigma})) \right) =B= M^0.$$
which comes from (\ref{*5}) and the following inequality (due to (\ref{*5}))
$$\min_{\sigma \in \{ +,- \}} \ (-H^{-\sigma}_{\alpha}(L^{\alpha,\sigma})) \geq \min_{\sigma \in \{ +,- \}} \ (-H^{-\sigma}_{\alpha}(p^{\alpha,0}_{i,\sigma})), \quad \mbox{for any} \quad i\geq 1.$$
\end{itemize}

{ \bf Step 3: Conclusion}\\
The estimates (\ref{spacetimegrad}) directly follow from (\ref{*10}), (\ref{*11}) and (\ref{*6}) and Theorem \ref{gradientestimate}.
\end{Proof}

\begin{Proof}{\bf of Theorem~\ref{convergence2}:}\\
{\bf Step 0: Preliminaries}\\
Let $T>0$ be fixed and let $\varepsilon:=(\Delta t, \Delta x)$ satisfy the CFL condition (\ref{CFL'}).

Assume that $u_\varepsilon$ is the numerical solution of (\ref{NS})-(\ref{CINS}).
We consider $\overline{u}$ and $\underline{u}$ respectively defined in (\ref{subsolution}) and (\ref{supersolution}).
By construction, we have
$$\underline{u} \leq \overline{u}.$$
We will show in the following steps that $\underline{u}$ (resp. $\overline{u}$) is a viscosity supersolution (resp. viscosity subsolution) of equation (\ref{HJ})-(\ref{CI}),
such that there exists a constant $C_T > 0$ such that for all $(t,y) \in [0,T) \times J$
$$\underline{u}(t,y) \geq -C_T(1+|y|) \quad \left( \mbox{resp.} \quad \overline{u}(t,y) \leq C_T(1+|y|) \right),$$
and such that
$$\underline{u}(0,y)\ge u_0(y) \quad \left( \mbox{resp.} \quad \overline{u}(0,y) \le u_0(y) \right) \quad \mbox{for all} \quad y \in J.$$
Using the comparison principle (Proposition~\ref{pro:1}), we obtain
$$\overline{u} \leq u \leq \underline{u}.$$
Thus from Definition~\ref{defviscosity}, we can conclude that $\overline{u}=u=\underline{u}$.
This implies the statement of Theorem~\ref{convergence2}.\\

{\bf Step 1: First bounds on the half relaxed limits}\\
From Lemma~\ref{lem:1}, we deduce that for any $\alpha=1,...,N$, any $i \geq0$ and any $n \geq 0$, we have
$$m^0_0 n \Delta t \leq U^{\alpha,n}_i -U^{\alpha,0}_i \leq M^0_0 n \Delta t.$$
Passing to the limit with $\varepsilon \to 0$ (always satisfying CFL condition (\ref{CFL'})), we get
$$u_0(y) + m^0_0 t \leq \underline{u}(t,y) \leq \overline{u}(t,y) \leq u_0(y) + M^0_0 t, \quad \mbox{for} \quad (t,y) \in [0,T) \times J.$$
This implies that
\begin{equation}\label{viscineq0}
\overline{u}(0,y) \leq u_0(y) \leq \underline{u}(0,y), \quad \mbox{for all} \quad y \in J,
\end{equation}
and
$$\overline{u}(t,y) \leq C_T(1+|y|) \quad \mbox{and} \quad \underline{u}(t,y) \geq -C_T(1+|y|),$$
with $C_T=\max \left\{|m^0_0|T,|M^0_0|T \right\} +|u_0(0)| +L$ and $L=\displaystyle \max_{\alpha, \sigma \in \{\pm\}} |L^{\alpha,\sigma}|$.

In next step, we show that $\underline{u}$ is a supersolution of (\ref{HJ})-(\ref{CI}) in the viscosity sense. We skip the proof that $\overline{u}$ is a viscosity subsolution because it is similar.\\

{\bf Step 2: Proof of $\underline{u}$ being a viscosity supersolution}\\
Let us consider $\underline{u}=\left(\underline{u}^{\alpha} \right)_{\alpha=1,...,N}$ as defined in (\ref{supersolution}) and a test function $\varphi:=(\varphi^\alpha)_{\alpha} \in \mathcal{C}^{1}_{*}([0,T) \times J)$ satisfying
$$ \begin{cases}
\underline{u} \geq \varphi \quad \mbox{on} \quad [0,T) \times J ,\\
\underline{u} = \varphi \quad \mbox{at} \quad P_0 = (t_0, y_0) \in [0,T) \times J, \quad \mbox{with} \quad t_0>0.
\end{cases} $$

Thus up to replace $\varphi(P)$ by $\hat{\varphi}(P) = \varphi(P) +|P-P_0|^2$, we can assume that
$$ \begin{cases}
\underline{u} > \varphi \quad \mbox{on} \quad \overline{B_{r}(P_0)} \setminus \{ P_0 \}, \\
\underline{u} = \varphi \quad \mbox{at} \quad P_0.
\end{cases} $$

We set $B_{r}(P_0)$ the ball centred at $P_0$ with fixed radius $r>0$, and set $\Omega_\varepsilon$ defined as the intersection between the closed ball centred on $P_0$ and the grid points $G_{\varepsilon}$ (defined in (\ref{grid})), i.e.
$$\Omega_{\varepsilon} := \overline{B_r(P_0)} \cap G_{\varepsilon}.$$
Note that for $\varepsilon$ small enough, we have $\Omega_{\varepsilon} \neq \emptyset$.
Up to decrease $r$, we can assume that $B_r(P_0) \subset [0,T-r) \times J$.

Define also
$$ M_{\varepsilon} = \displaystyle{\inf_{\Omega_{\varepsilon}}} \ (u_{\varepsilon} - {\varphi}) = (u_{\varepsilon} - {\varphi})(P_{\varepsilon}) ,$$
where
$$P_{\varepsilon}=(t_{\varepsilon},y_{\varepsilon}) \in [0,T)\times J_{\alpha_{\varepsilon}} \quad \text{with} \quad y_{\varepsilon}=x_{\varepsilon}e_{\alpha_{\varepsilon}} \quad \text{and} \quad \begin{cases}
t_{\varepsilon}:=n_{\varepsilon} \Delta t \\
x_{\varepsilon}:=i_{\varepsilon} \Delta x
\end{cases}.$$

By the definition of $\underline{u}$ in (\ref{supersolution}), it is classical to show that if $\varepsilon \to 0$ we get the following (at least for a subsequence)
\begin{equation}\label{lim}
\begin{cases}
M_{\varepsilon}=(u_{\varepsilon}-\varphi)(P_{\varepsilon}) \to M_0 = \displaystyle{\inf_{\overline{B_r(P_0)}}} \ (\underline{u}-{\varphi})=0 ,\\
P_{\varepsilon} \to P_0.
\end{cases}
\end{equation}


Let us now check that $\underline{u}$ is a viscosity supersolution of (\ref{HJ}).
To this end, using Proposition~\ref{pro:4} we want to show that
\begin{itemize}
\item[$\bullet$] if $y_0=x_0 e_{\alpha_0} \in J_{\alpha_0}^{*}$ for a given $\alpha_0$, then
$$ \varphi^{\alpha_0}_t+H_{\alpha_0}(\varphi^{\alpha_0}_x) \geq 0 \quad \mbox{at} \quad (t_0,x_0),$$
\item[$\bullet$] if $y_0=0$, then either there exists one index $\alpha_0$ such that
$$ \varphi_t^{\alpha_0}+H_{\alpha_0}(\varphi_x^{\alpha_0}) \geq 0 \quad \mbox{at} \quad (t_0,0),$$
or we have $$\varphi_t + \displaystyle{\sup_{\alpha=1,...,N}} \ \left\{ H_{\alpha}^{-}(\varphi^{\alpha}_x) \right\} \geq 0 \quad \mbox{at} \quad (t_0,0).$$
\end{itemize}

Because $t_0>0$ and $P_{\varepsilon} \to P_0$, this implies in particular that $n_{\varepsilon} \geq 1$ for $\varepsilon$ small enough.
We have to distinguish two cases according to the value of $y_0$.

{\bf Case 1: $P_0=(t_0, y_0)$ with $y_0=0$}

We distinguish two subcases, up to subsequences.

{\small {\bf Subcase 1.1: $P_{\varepsilon}=(t_{\varepsilon}, y_{\varepsilon})$ with $y_{\varepsilon}=y_0=0$}}

Using the definitions (\ref{defSN1}), (\ref{defSN2}) and the numerical scheme (\ref{NS}), we recall that for all $n \geq 0$ and for any $\alpha=1,...,N$
$$\begin{aligned}
U^{\alpha,n+1}_{0}=:U^{n+1}_{0}&=U^{n}_{0} - \Delta t \displaystyle{\max_{\alpha=1,\ldots ,N}} \ H^{-}_{\alpha} \left(\dfrac{U^{\alpha,n}_{1}-U^{n}_{0}}{\Delta x}\right) \\
&= S_0 \left[U^{n}_{0}, \left(U^{\alpha,n}_{1} \right)_{\alpha=1,...,N}\right]
\end{aligned}$$
where $S_0$ is monotone under the CFL condition (\ref{CFL'}) (see Proposition~\ref{pro1}).

Let ${\varphi}_{\varepsilon}:=M_{\varepsilon}+{\varphi}$ such that
$$ \begin{aligned}
{\varphi}_{\varepsilon}(P_{\varepsilon}) =u_{\varepsilon}(P_{\varepsilon}) &= U^{n_{\varepsilon}}_{0} \\
&= S_0 \left[U^{n_{\varepsilon}-1}_0,\left(U^{\alpha,n_{\varepsilon}-1}_1 \right)_{\alpha=1,...,N}\right] \\
&\geq S_0 \left[{\varphi}_{\varepsilon} ( (n_{\varepsilon}-1) \Delta t , 0), \left({\varphi}^{\alpha}_{\varepsilon} ( (n_{\varepsilon}-1) \Delta t , \Delta x)\right)_{\alpha=1,...,N} \right]
\end{aligned} $$
where we use the monotonicity of the scheme in the last line and the fact that $u_{\varepsilon} \geq \varphi_{\varepsilon}$ on $\Omega_{\varepsilon}$.

Thus we have
$$ \dfrac{{\varphi}_{\varepsilon} ( n_{\varepsilon}\Delta t , 0) - {\varphi}_{\varepsilon} ( (n_{\varepsilon}-1) \Delta t , 0)}{\Delta t} + \displaystyle{\max_{\alpha=1,...,N}} H^{-}_{\alpha}\left(\dfrac{\varphi^{\alpha}_{\varepsilon} ( (n_{\varepsilon}-1) \Delta t , \Delta x)- {\varphi}^{\alpha}_{\varepsilon} ( (n_{\varepsilon}-1) \Delta t , 0)}{\Delta x}\right) \geq 0 .$$

This implies
$$ (\varphi_{\varepsilon})_t + \displaystyle{\max_{\alpha=1,\ldots,N}} \ H^{-}_{\alpha} ((\varphi^{\alpha}_{\varepsilon})_x) + o_{\varepsilon}(1) \geq 0 \quad \mbox{at} \quad (t_{\varepsilon},0)$$

and passing to the limit with $\varepsilon \to 0$, we get the supersolution condition at the junction point
\begin{equation}\label{*1}
\varphi_t + \displaystyle{\max_{\alpha=1,\ldots,N}} \ H^{-}_{\alpha} (\varphi^{\alpha}_x) \geq 0 \quad \mbox{at} \quad (t_0,0).
\end{equation}

{\small {\bf Subcase 1.2: $P_{\varepsilon}=(t_{\varepsilon}, y_{\varepsilon})$ with $y_{\varepsilon} \in J^{*}_{\alpha_{\varepsilon}}$}}

In this case, the infimum $M_{\varepsilon}$ is reached for a point on the branch $\alpha_{\varepsilon}$ which is different from the junction point. Thus the definitions (\ref{defSN1}), (\ref{defSN2}) and the numerical scheme (\ref{NS}) give us that for all $n \geq 0$ and $i \geq 1$
$$\begin{aligned}
U^{\alpha_{\varepsilon},n+1}_{i} &= U^{\alpha_{\varepsilon},n}_{i} + \Delta t \min \ \{-H^{-}_{\alpha_{\varepsilon}} (p^{\alpha_{\varepsilon},n}_{i,+}), -H^{+}_{\alpha_{\varepsilon}} (p^{\alpha_{\varepsilon},n}_{i,-}) \} \\
&= S_{\alpha_{\varepsilon}} [U^{\alpha_{\varepsilon},n}_{i-1},U^{\alpha_{\varepsilon},n}_{i},U^{\alpha_{\varepsilon},n}_{i+1}].
\end{aligned}$$

Let ${\varphi}^{\alpha_{\varepsilon}}_{\varepsilon}:=M_{\varepsilon}+{\varphi}^{\alpha_{\varepsilon}}$ such that 
$$ \begin{aligned}
{\varphi}^{\alpha_{\varepsilon}}_{\varepsilon}(P_{\varepsilon})=u^{\alpha_\varepsilon}_{\varepsilon}(P_{\varepsilon}) &= U^{\alpha_{\varepsilon},n_{\varepsilon}}_{i_{\varepsilon}} \\
&= S_{\alpha_{\varepsilon}} [U^{\alpha_{\varepsilon},n_{\varepsilon}-1}_{i_{\varepsilon}-1}, U^{\alpha_{\varepsilon},n_{\varepsilon}-1}_{i_{\varepsilon}}, U^{\alpha_{\varepsilon},n_{\varepsilon}-1}_{i_{\varepsilon}+1}] \\
&\geq S_{\alpha_{\varepsilon}} [{\varphi}^{\alpha_{\varepsilon}}_{\varepsilon} ( (n_{\varepsilon}-1) \Delta t , (i_{\varepsilon}-1) \Delta x), {\varphi}^{\alpha_{\varepsilon}}_{\varepsilon} ( (n_{\varepsilon}-1) \Delta t , i_{\varepsilon} \Delta x), {\varphi}^{\alpha_{\varepsilon}}_{\varepsilon} ( (n_{\varepsilon}-1) \Delta t , (i_{\varepsilon}+1) \Delta x) ]
\end{aligned} $$
where we use the monotonicity of the scheme and the fact that $u^{\alpha_{\varepsilon}}_{\varepsilon} \geq \varphi^{\alpha_{\varepsilon}}_{\varepsilon}$ in the neighbourhood of $(t_{\varepsilon},x_{\varepsilon})$.

Thus we have that for any $\varepsilon=(\Delta t, \Delta x)$
$$\begin{aligned}
0 \leq & \dfrac{{\varphi}^{\alpha_{\varepsilon}}_{\varepsilon} (n_{\varepsilon}\Delta t , i_{\varepsilon} \Delta x) - {\varphi}^{\alpha_{\varepsilon}}_{\varepsilon} ((n_{\varepsilon}-1) \Delta t , i_{\varepsilon} \Delta x)}{\Delta t} \\
&+ \max \Bigg\{ H^+_{\alpha_{\varepsilon}}\left(\dfrac{\varphi^{\alpha_{\varepsilon}}_{\varepsilon} (n_{\varepsilon}\Delta t,i_{\varepsilon} \Delta x)- {\varphi}^{\alpha_{\varepsilon}}_{\varepsilon} (n_{\varepsilon}\Delta t, (i_{\varepsilon}-1) \Delta x)}{\Delta x}\right),\\
& \qquad \qquad H^-_{\alpha_{\varepsilon}}\left(\dfrac{\varphi^{\alpha_{\varepsilon}}_{\varepsilon} (n_{\varepsilon}\Delta t,(i_{\varepsilon}+1) \Delta x)- {\varphi}^{\alpha_{\varepsilon}}_{\varepsilon} (n_{\varepsilon}\Delta t, i_{\varepsilon} \Delta x)}{\Delta x}\right) \Bigg\} .
\end{aligned}$$

Since $H_{\alpha}(p)=\max \left\{ H^+_{\alpha}(p),H^-_{\alpha}(p) \right\}$, this implies
$$ (\varphi^{\alpha_{\varepsilon}})_t + H_{\alpha_{\varepsilon}} ((\varphi^{\alpha_{\varepsilon}})_x) + o_{\varepsilon}(1) \geq 0 \quad \mbox{at} \quad (t_{\varepsilon},x_{\varepsilon}).$$

Up to a subsequence, we can assume that $\alpha_{\varepsilon}$ is independent of $\varepsilon$ and equal to $\alpha_0$. Thus passing to the limit with $\varepsilon \to 0$, we obtain
\begin{equation}\label{*2}
{\varphi}^{\alpha_0}_t + H_{\alpha_0} ({\varphi}^{\alpha_0}_x) \geq 0 \quad \mbox{at} \quad (t_0,0).
\end{equation}

{ \bf Case 2: $P_0=(t_0,y_0)$ with $y_0=x_0 e_{\alpha_0} \in J^{*}_{\alpha_0}$}

As $P_{\varepsilon} \to P_0$ from (\ref{lim}), we can always consider that for $\varepsilon$ small enough, we can write $P_{\varepsilon}=(t_{\varepsilon}, y_{\varepsilon})$ with $y_{\varepsilon} \in J^{*}_{\alpha_{\varepsilon}}$.
Thus the proof for this case is similar to the one in Subcase 1.2.
We then conclude
\begin{equation}\label{*3}
{\varphi}^{\alpha_0}_t + H_{\alpha_0} ({\varphi}^{\alpha_0}_x) \geq 0 \quad \mbox{at} \quad (t_0,x_0).
\end{equation}

{\bf Step 3: Conclusion}\\
The results (\ref{viscineq0}), (\ref{*1}), (\ref{*2}) and (\ref{*3}) imply that $\underline{u}$ is a viscosity supersolution of (\ref{HJ})-(\ref{CI}).
This ends the proof of Theorem~\ref{convergence2}.
\end{Proof}

%

\begin{Proof}{\bf of Theorem~\ref{convergence1}:}
The proof is quite similar to the proof of Theorem~\ref{convergence2}. However it differs on some points mainly because we do not know if the comparison principle from Proposition~\ref{pro:1} holds for (\ref{HJ}).
\begin{itemize}
\item[$\bullet$] We recall from Lemma~\ref{lem:1} that $u^{\alpha}_{\varepsilon}(n \Delta t, i \Delta x):=U^{\alpha,n}_i$ with $\varepsilon=(\Delta t, \Delta x)$ enjoys some discrete Lipschitz bounds in time and space, independent of $\varepsilon$.
\item[$\bullet$] It is then possible to extend the discrete function $u_{\varepsilon}$, defined only on the grid points, into a continuous function $\tilde{u}_{\varepsilon}$, with the $\mathcal{Q}_1$ quadrilateral finite elements approximation for which we have the same Lipschitz bounds. We recall that the approximation is the following: consider a map $(t,x)\mapsto u(t,x)$ that takes values only on the vertex of a rectangle $ABCD$ with $A=(0,0)$, $B=(0,1)$, $C=(1,1)$ and $D=(1,0)$ (for sake of simplicity we take $\Delta t=1=\Delta x$). Then we extend the map $u$ to any point $(t,x)$ of the rectangle such that
$$u(t,x)=[u_A+x(u_B-u_A)](1-t)+[u_D+x(u_C-u_D)]t.$$
\item[$\bullet$] In this way we can apply the Ascoli theorem which shows that there exists a subsequence $\tilde{u}_{\varepsilon'}$ which converges towards a function $u$, uniformly on every compact set (in time and space).
\item[$\bullet$] We can then conclude that $u$ is a viscosity super and subsolution of (\ref{HJ})-(\ref{CI}) repeating the proof of Theorem~\ref{convergence2}.
\end{itemize}
This ends the proof.
\end{Proof}

\begin{Proof}{\bf of Corollary \ref{cor2}:}
The proof combines the gradient and time derivative estimates from Lemma~\ref{lem:1} and the results of convergence from Theorem~\ref{convergence1}. Indeed, passing to the limit in (\ref{spacetimegrad}) for a subsequence $\varepsilon'$, using the convergence result of Theorem~\ref{convergence1}, we finally get (\ref{contestimates}).
\end{Proof}

\section{Application to traffic flow}
\label{sect::application}

As our motivation comes from traffic flow modelling, this section is devoted to the traffic interpretation of the model and the scheme. Notice that \cite{IMZ} has already focused on the meaning of the junction condition in this framework.

\subsection{Settings}
We first recall the main variables adapted for road traffic modelling as they are already defined in \cite{IMZ}. 
We consider a junction with $N_I \geq 1$ incoming roads and $N_O \geq 1$ outgoing ones. 
We also set that $N_I+N_O=:N$.\\

\textbf{Densities and scalar conservation law.} \ We assume that the vehicles densities denoted by $(\rho^{\alpha})_{\alpha}$ solve the following scalar conservation laws (also called LWR model for Lighthill, Whitham \cite{LW} and Richards \cite{Richards}):
\begin{equation}\label{LWRnet}
\begin{cases}
\rho_t^{\alpha} + (f^{\alpha}(\rho^{\alpha}))_X=0, \quad \mbox{for} \quad (t,X) \in [0,+\infty) \times (-\infty,0), \quad {\alpha}=1,...,N_I,\\
\rho_t^{\alpha} + (f^{\alpha}(\rho^{\alpha}))_X=0, \quad \mbox{for} \quad (t,X) \in [0,+\infty) \times (0,+\infty), \quad {\alpha}=N_I+1,...,N_I+N_O,
\end{cases}
\end{equation}
where we assume that the junction point is located at the origin $X=0$.

We assume that for any $\alpha$ the flux function $f^\alpha : \R \to \R$ reaches its unique maximum value for a critical density $\rho=\rho_{c}^\alpha >0$ and it is non decreasing on $(-\infty,\rho^{\alpha}_c)$ and non-increasing on $(\rho^{\alpha}_c,+\infty)$. In traffic modelling, $\rho^\alpha \mapsto f^\alpha(\rho^\alpha)$ is usually called the \textit{fundamental diagram}.

Let us define for any $\alpha=1,...,N$ the Demand function $f^\alpha_D$ (resp. the Supply function $f^\alpha_S$) such that
\begin{equation}\label{defdemandsupply}
f^\alpha_D(p)=
\begin{cases}
f^\alpha(p) \quad &\mbox{for} \quad p \leq \rho^\alpha_c \\
f^\alpha(\rho^\alpha_c) \quad &\mbox{for} \quad p \geq \rho^\alpha_c
\end{cases}
\qquad \left( \mbox{resp.} \quad
f^\alpha_S(p)=
\begin{cases}
f^\alpha(\rho^\alpha_c) \quad &\mbox{for} \quad p \leq \rho^\alpha_c \\
f^\alpha(p) \quad &\mbox{for} \quad p \geq \rho^\alpha_c
\end{cases} \right).
\end{equation}

 We assume that we have a set of fixed coefficients $0 \leq \left( \gamma^{\alpha} \right)_{\alpha} \leq 1$ that denote:
\begin{itemize}
\item[$\bullet$] either the proportion of the flow from the branch $\alpha=1,...,N_I$ which enters in the junction,
\item[$\bullet$] or the proportion of the flow on the branch $\alpha=N_I+1,...,N$ exiting from the junction.
\end{itemize}
We also assume the natural relations
$$\sum_{\alpha=1}^{N_I} \gamma^{\alpha} =1 \quad \mbox{and} \quad \sum_{\beta=N_I+1}^{N_I+N_O} \gamma^{\beta} =1.$$

\begin{rem}
We consider that the coefficients $(\gamma^\alpha)_{\alpha=1,...,N}$ are fixed and known at the beginning of the simulations. Such framework is particularly relevant for ``quasi stationary'' traffic flows.
\end{rem}

\textbf{Vehicles labels and Hamilton-Jacobi equations.} \ Extending for any $N_I \geq 1$ the interpretation and the notations given in \cite{IMZ} for a single incoming road, let us consider the \textit{continuous} analogue $u^\alpha$ of the discrete vehicles labels (in the present paper with labels increasing in the backward direction with respect to the flow)
\begin{equation}\label{defu}
\left\{ \begin{array}{lll}
u^{\alpha} (t,x) = u(t,0) - \dfrac{1}{\gamma^{\alpha}} \displaystyle{\int_0^{-x}} \rho^{\alpha}(t,y) dy, \quad &\mbox{for} \quad x:=-X>0, \quad &\mbox{if} \quad \alpha=1,...,N_I ,\\
\\
u^{\beta} (t,x) = u(t,0) - \dfrac{1}{\gamma^{\beta}} \displaystyle{\int_0^x} \rho^{\beta}(t,y) dy, \quad &\mbox{for} \quad x:=X>0, \quad &\mbox{if} \quad \beta=N_I+1,...,N ,
\end{array} \right.
\end{equation}
with equality of the functions at the junction point ($x=0$), i.e.
\begin{equation}\label{eqatjunction}
u^{\alpha}(t,0)=u^{\beta}(t,0)=:u(t,0) \quad \mbox{for any} \quad \alpha, \beta .
\end{equation}
where the common value $u(t,0)$ is nothing else than the (continuous) label of the vehicle at the junction point.

Following \cite{IMZ}, for a suitable choice of the function $u(t,0)$, it is possible to check that the vehicles labels $u^{\alpha}$ satisfy the following Hamilton-Jacobi equation:
\begin{equation}\label{HJ*}
u^{\alpha}_t + H_{\alpha}(u^{\alpha}_x) = 0, \quad \mbox{for} \quad (t,x) \in [0,+\infty) \times (0,+\infty), \quad {\alpha}=1,...,N
\end{equation}
where
\begin{equation}\label{Hfunctf}
H_{\alpha}(p) := \begin{cases} 
-\dfrac{1}{\gamma^{\alpha}} f^{\alpha}(\gamma^{\alpha}p) \quad &\mbox{for} \quad \alpha=1,...,N_I ,\\
\\
-\dfrac{1}{\gamma^{\alpha}} f^{\alpha}(-\gamma^{\alpha}p) \quad &\mbox{for} \quad \alpha=N_I+1,...,N_I+N_O.
\end{cases}
\end{equation}

Following definitions of $H^-_{\alpha}$ and $H^+_{\alpha}$ in  (\ref{Hamiltonians+-}) we get
\begin{equation}\label{defsupplydemand}
H^-_{\alpha}(p)=
\begin{cases}
-\dfrac{1}{\gamma^\alpha}f^\alpha_D(\gamma^{\alpha}p) \quad &\mbox{for} \quad \alpha \leq N_I,\\
\\
-\dfrac{1}{\gamma^\alpha}f^\alpha_S(-\gamma^{\alpha}p) \quad &\mbox{for} \quad \alpha \geq N_I+1,
\end{cases} \quad \mbox{and} \quad 
H^+_{\alpha}(p)=
\begin{cases}
-\dfrac{1}{\gamma^\alpha}f^\alpha_S(\gamma^{\alpha}p) \quad &\mbox{for} \quad \alpha \leq N_I,\\
\\
-\dfrac{1}{\gamma^\alpha}f^\alpha_D(-\gamma^{\alpha}p) \quad &\mbox{for} \quad \alpha \geq N_I+1.
\end{cases}
\end{equation}

The junction condition in (\ref{HJ}) that reads
\begin{equation}\label{junction cond}
u_t (t,0) + \displaystyle{\max_{\alpha=1,...,N} H_{\alpha}^{-}(u_x(t,0^+))} =0.
\end{equation}
is a natural condition from the traffic point of view.
Indeed condition (\ref{junction cond}) can be rewritten as
\begin{equation}\label{PR4}
\begin{aligned}
u_t (t,0) = \min \left( \displaystyle{\min_{\alpha=1,...,N_I} \dfrac{1}{\gamma^{\alpha}} f^{\alpha}_D(\rho^{\alpha} (t,0^-))} , \displaystyle{\min_{\beta=N_I+1,...,N} \dfrac{1}{\gamma^{\beta}} f^{\beta}_S(\rho^{\beta} (t,0^+))} \right).
\end{aligned}
\end{equation}

The condition (\ref{PR4}) claims that the passing flux is equal to the minimum between the upstream demand and the downstream supply functions as they were presented by Lebacque in \cite{Lebacque93} and \cite{Lebacque96} (also for the case of junctions). This condition maximises the flow through the junction. This is also related to the Riemann solver \textit{RS2} in \cite{GP09} for junctions.

\subsection{Review of the literature with application to traffic}

\textbf{Junction modelling.} \ There is an important and fast growing literature about junction modeling from a traffic engineering viewpoint: see \cite{KL,TCCI,FR} for a critical review of junction models. The literature mainly refers to \textit{pointwise junction} models \cite{KL,LK02,LK05}. Pointwise models are commonly restated in many instances as optimization problems.\\

{\bf Scalar one dimensional conservation laws and networks.} \ Classically, macroscopic traffic flow models are based on a scalar one dimensional conservation law, e.g. the so-called LWR model (Lighthill, Whitham \cite{LW} and Richards \cite{Richards}). The literature is also quite important concerning hyperbolic systems of conversation laws (see for example \cite{Bressan,Dafermos,Lax,Serre} and references therein) but these books also propose a large description of the scalar case. It is well-known that under suitable assumptions there exists a unique weak entropy solution for scalar conservation laws without junction. \\

Until now, existence of weak entropy solutions for a Cauchy problem on a network has been proved for general junctions in \cite{GP09}. See also Garavello and Piccoli's book \cite{GP}. Uniqueness for scalar conservation laws for a junction with two branches has been proved first in \cite{GNPT} and then in \cite{AKR} under suitable assumptions. Indeed \cite{AKR} introduces a general framework with the notion of \textit{$L^1$-dissipative admissibility germ} that is a selected family of elementary solutions. To the best authors' knowledge, there is no uniqueness result for general junctions.\\

The conservation law counterpart of model (\ref{HJ*}),(\ref{eqatjunction}),(\ref{junction cond}) has been studied in \cite{GP09} as a Riemann solver called \textit{RS2}. In \cite{GP09} an existence result is presented for concave flux functions, using the Wave Front Tracking (WFT) method. Moreover the Lipschitz continuous dependence of the solution with respect to the initial data is proven. This shows that the process of construction of a solution (here the WFT method) creates a single solution. Nevertheless, up to our knowledge, there is no differential characterisation of this solution. Therefore the uniqueness of this solution is still an open problem.\\

{\bf Numerical schemes for conservation laws.} \ According to \cite{GR} and \cite{Leveque}, the numerical treatment of scalar conservation laws mainly deals with first order numerical schemes based on upwind finite difference method, such as the Godunov scheme \cite{Godunov} which is well-adapted for the LWR model \cite{Lebacque96}.

As finite difference methods introduce numerical viscosity, other techniques were developed such as kinetic schemes that derive from the kinetic formulation of hyperbolic equations \cite{Perthame}. Such kinetic schemes are presented in \cite{ADN} and they are applied to the traffic case in \cite{BNP06,BNP07,BNP07b}.\\

In \cite{GHZ} the authors apply a semidiscrete central numerical scheme to the Hamilton-Jacobi formulation of the LWR model. The equivalent scheme for densities recovers the classical Lax-Friedrichs scheme. Notice that the authors need to introduce at least two \textit{ghost-cells} on each branch near the junction point to counterstrike the dispersion effects when computing the densities at the junction.

\subsection{Derived scheme for the densities}
\label{NScons}

The aim of this subsection is to properly express the numerical scheme satisfied by the densities in the traffic modelling framework. Let us recall that the density denoted by $\rho^\alpha$ is a solution of (\ref{LWRnet}). 

\begin{figure}[!ht]
\begin{center}
\resizebox{5cm}{!}{\input{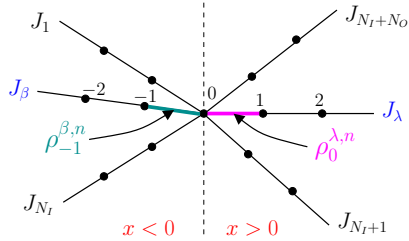}}
\caption{Discretization of the branches with the nodes for $(U^{\alpha,n}_i)$ and the segments for $(\rho^{\alpha,n}_i)$}
\label{fig:0}
\end{center}
\end{figure}

Let us consider a discretization with time step $\Delta t$ and space step $\Delta x$. Then we define the discrete car density $\rho^{\alpha,n}_i \geq 0$ for $n \geq 0$ and $i \in \Z$ (see Figure \ref{fig:0}) by
\begin{equation}\label{defrho}
\rho^{\alpha,n}_i := \left\{
\begin{array}{lll}
\gamma^{\alpha} p^{\alpha,n}_{|i|-1,+}, \quad &\mbox{for} \quad i \leq -1, \quad &\alpha=1,...,N_I, \\ 
\\
-\gamma^{\alpha} p^{\alpha,n}_{i,+}, \quad &\mbox{for} \quad i \geq 0, \quad &\alpha=N_I+1,...,N_I+N_O ,
\end{array} \right.
\end{equation}
where we recall
$$p^{\alpha,n}_{j,+}:= \dfrac{U^{\alpha,n}_{j+1}-U^{\alpha,n}_j}{\Delta x}, \quad \mbox{for} \quad j \in \N, \quad \alpha=1,...,N.$$

We have the following result:
\begin{lem}[Derived numerical scheme for the density]
\label{derivsch}
\ \\
If $(U^{\alpha,n}_i)$ stands for the solution of (\ref{NS})-(\ref{CINS}), then the density $(\rho^{\alpha,n}_i)$ defined in (\ref{defrho}) is a solution of the following numerical scheme for $\alpha=1,...,N$
 \begin{equation}\label{LWRscheme}
\dfrac{\Delta x}{\Delta t}\{ \rho^{\alpha,n+1}_{i}-\rho^{\alpha,n}_{i} \}=
\begin{cases}
F^\alpha(\rho^{\alpha,n}_{i-1},\rho^{\alpha,n}_{i}) - F^\alpha(\rho^{\alpha,n}_{i},\rho^{\alpha,n}_{i+1}) \quad &\mbox{for} \ \begin{cases}
i \leq -1 \ &\mbox{if} \quad \alpha \leq N_I, \\
i \geq 1 \ &\mbox{if} \quad \alpha \geq N_I+1,
\end{cases} \\
\\
F^\alpha_0\left((\rho^{\beta,n}_{-1})_{\beta \leq N_I},(\rho^{\lambda,n}_{0})_{\lambda \geq N_I+1}\right) - F^\alpha(\rho^{\alpha,n}_{i},\rho^{\alpha,n}_{i+1}) \quad &\mbox{for} \quad i = 0, \quad \alpha \geq N_I+1,\\
\\
F^\alpha(\rho^{\alpha,n}_{i-1},\rho^{\alpha,n}_{i}) - F^\alpha_0\left((\rho^{\beta,n}_{-1})_{\beta \leq N_I},(\rho^{\lambda,n}_{0})_{\lambda \geq N_I+1}\right) \quad &\mbox{for} \quad i = -1, \quad \alpha \leq N_I,
\end{cases}
 \end{equation}

where we define the fluxes by
\begin{equation}\label{flux}
\begin{cases}
 F^\alpha(\rho^{\alpha,n}_{i-1},\rho^{\alpha,n}_{i}) := \min \left\{f_D^\alpha(\rho^{\alpha,n}_{i-1}), \ f_S^\alpha(\rho^{\alpha,n}_{i}) \right\} \quad &\mbox{for} \quad
\begin{cases}
i \leq -1 \quad &\mbox{if} \quad \alpha \leq N_I, \\
i \geq 1 \quad &\mbox{if} \quad \alpha \geq N_I+1,
\end{cases} \\
 \\
 F^\alpha_0\left((\rho^{\beta,n}_{-1})_{\beta \leq N_I},(\rho^{\lambda,n}_{0})_{\lambda \geq N_I+1}\right) := \gamma^{\alpha} F_0 \quad &\mbox{for} \quad \alpha=1,...,N, \\
 \\
 F_0 := \min \left\{ \displaystyle{\min_{\beta\leq N_I}} \ \dfrac{1}{\gamma^\beta} f_D^\beta(\rho^{\beta,n}_{-1}), \ \displaystyle{\min_{\lambda \geq N_I+1}} \ \dfrac{1}{\gamma^\lambda} f_S^\lambda(\rho^{\lambda,n}_{0}) \right\}.
 \end{cases}
 \end{equation}
and $f_S^\alpha$, $f_D^\alpha$ are defined in (\ref{defdemandsupply}).

The initial condition is given by
\begin{equation}\label{CI-LWR}
\rho^{\alpha,0}_i := \left\{
\begin{array}{lll}
\gamma^{\alpha} \dfrac{u^{\alpha}_{0}(|i| \Delta x) -u^{\alpha}_{0}((|i|-1) \Delta x)}{\Delta x} , \quad &\mbox{for} \quad i \leq -1, \quad &\alpha=1,...,N_I, \\ 
\\
\gamma^{\alpha} \dfrac{u^{\alpha}_{0}(i \Delta x) -u^{\alpha}_{0}((i+1) \Delta x)}{\Delta x} , \quad &\mbox{for} \quad i \geq 0, \quad &\alpha=N_I+1,...,N_I+N_O.
\end{array} \right.
\end{equation}
\end{lem}

\begin{rem}
Notice that (\ref{LWRscheme}) recovers the classical Godunov scheme \cite{Godunov} for $i \neq 0, -1$ while it is not standard for the two other cases $i=0, -1$. Moreover we can check that independently of the chosen CFL condition, the scheme (\ref{LWRscheme}) is not monotone (at the junction, $i=0$ or $i=-1$) if the total number of branches $N \geq 3$ and is monotone if $N=2$ for a suitable CFL condition. 
\end{rem}

\begin{rem}
From (\ref{defcontgrad}), (\ref{NS}) and (\ref{defsupplydemand}), we can easily show that
$$\begin{aligned}
m^0 = \min \Bigg\{ &\min_{\substack{\alpha \leq N_I \\ i \leq -1}} \min \left( \dfrac{1}{\gamma^\alpha} f^{\alpha}_D(\rho^{\alpha,0}_{i-1}), \dfrac{1}{\gamma^\alpha} f^{\alpha}_S(\rho^{\alpha,0}_{i}) \right), \\
&\min_{\substack{\alpha \geq N_I+1 \\ i \geq 1}} \min \left( \dfrac{1}{\gamma^\alpha} f^{\alpha}_D(\rho^{\alpha,0}_{i-1}), \dfrac{1}{\gamma^\alpha} f^{\alpha}_S(\rho^{\alpha,0}_{i}) \right), \\
&\min \left( \min_{\alpha \leq N_I} \dfrac{1}{\gamma^\alpha} f^{\alpha}_D(\rho^{\alpha,0}_{-1}), \min_{\alpha \geq N_I+1} \dfrac{1}{\gamma^\alpha} f^{\alpha}_S(\rho^{\alpha,0}_{0}) \right) \Bigg\},
\end{aligned}$$
with the first part dealing with incoming branches, the second with outgoing branches and the third with the junction point.
As $\displaystyle f^{\alpha}(p)=\min \left\{ f^{\alpha}_S(p), f^{\alpha}_D(p) \right\}$ for any $p$, the latter can be rewritten as the minimal initial flux
$$m^0 = \min \left\{ \min_{\substack{\alpha \leq N_I \\ i \leq -1}} \left( \dfrac{1}{\gamma^\alpha} f^{\alpha}(\rho^{\alpha,0}_{i}) \right), \min_{\substack{\alpha \geq N_I+1 \\ i \geq 0}} \left( \dfrac{1}{\gamma^\alpha} f^{\alpha}(\rho^{\alpha,0}_{i}) \right) \right\}.$$

We set for any $\alpha=1,...,N$
$$\begin{cases}
\rho^-_{\alpha} := \left( f^{\alpha}_D \right)^{-1} (\gamma^\alpha m^0),\\
\rho^+_{\alpha} := \left( f^{\alpha}_S \right)^{-1} (\gamma^\alpha m^0),
\end{cases}$$
From Theorem~\ref{gradientestimate} and Remark \ref{rem:2}, if (\ref{CFL'}) is satisfied then it is easy to check that
$$\rho^-_{\alpha} \leq \rho^{\alpha,n}_{i} \leq \rho^+_{\alpha}, \quad \mbox{for all} \quad n \geq 0.$$
Then the CFL condition (\ref{CFL'}) can be rewritten for the densities as
\begin{equation}\label{CFLrho}
\dfrac{\Delta x}{\Delta t} \geq \sup_{\substack{\alpha=1,...,N \\ \rho^{\alpha} \in [\rho^-_{\alpha},\rho^+_{\alpha}]}} \ \left|(f^{\alpha})'(\rho^{\alpha})\right|.
\end{equation}
\end{rem}

\begin{Proof}{\bf of Lemma~\ref{derivsch}:}
We distinguish two cases according to if we are either on an incoming or an outgoing branch. We investigate the incoming case. The outgoing case can be done similarly.

Let us consider any $\alpha=1,...,N_I$, $n\geq 0$ and $i \leq -1$.

According to (\ref{defrho}),  for $i \leq -2$ we have that:
$$\begin{aligned}
\dfrac{\rho^{\alpha,n+1}_{i}-\rho^{\alpha,n}_{i}}{\Delta t} &= \dfrac{\gamma^{\alpha}}{\Delta x \Delta t} \left\{ \left(U^{\alpha,n+1}_{|i|}-U^{\alpha,n+1}_{|i|-1}\right) - \left(U^{\alpha,n}_{|i|}-U^{\alpha,n}_{|i|-1}) \right) \right\}\\
&= \dfrac{\gamma^{\alpha}}{\Delta x} \left\{ \min \left(-H^-_{\alpha}(p^{\alpha,n}_{|i|,+}),-H^+_{\alpha}(p^{\alpha,n}_{|i|,-}) \right) - \min \left(-H^-_{\alpha}(p^{\alpha,n}_{|i|-1,+}),-H^+_{\alpha}(p^{\alpha,n}_{|i|-1,-}) \right) \right\} \\
&= \dfrac{1}{\Delta x} \left\{ \min \left(f^{\alpha}_D(\rho^ {\alpha,n}_{i-1}),f^{\alpha}_S(\rho^ {\alpha,n}_{i}) \right) - \min \left(f^{\alpha}_D(\rho^ {\alpha,n}_{i}),f^{\alpha}_S(\rho^ {\alpha,n}_{i+1}) \right) \right\}
\end{aligned}$$
where we use the numerical scheme (\ref{NS}) in the second line and (\ref{defsupplydemand}) in the last line.

We then recover the result if we set the fluxes functions $F^{\alpha}$ as defined in (\ref{flux}).\\

For the special case of $i=-1$, we have
$$\begin{aligned}
\dfrac{\rho^{\alpha,n+1}_{-1}-\rho^{\alpha,n}_{-1}}{\Delta t} &= \dfrac{\gamma^{\alpha}}{\Delta x} \left\{ \left(\dfrac{U^{\alpha,n+1}_{1}-U^{\alpha,n}_{1}}{\Delta t} \right) - \left( \dfrac{U^{\alpha,n+1}_{0}-U^{\alpha,n}_{0}}{\Delta t} \right) \right\} \\
&= \dfrac{\gamma^{\alpha}}{\Delta x} \left\{\min \left(-H^-_{\alpha}(p^{\alpha,n}_{1,+}),-H^+_{\alpha}(p^{\alpha,n}_{1,-}) \right)  - \min_{\beta=1,...,N} \left(-H^-_{\beta}(p^{\beta,n}_{0,+}) \right) \right\} \\
&= \dfrac{1}{\Delta x} \left\{ \min \left(f^{\alpha}_D(\rho^ {\alpha,n}_{-2}),f^{\alpha}_S(\rho^ {\alpha,n}_{-1}) \right) - \gamma^\alpha \min \left(\min_{\beta=1,...,N_I} \dfrac{1}{\gamma^\beta} f^{\beta}_D(\rho^ {\beta,n}_{-1}),\min_{\lambda=N_I+1,...,N} \dfrac{1}{\gamma^\lambda} f^{\lambda}_S(\rho^ {\lambda,n}_{0}) \right) \right\}
\end{aligned}$$
where in the last line we have used (\ref{defsupplydemand}). Setting the flux function $F^{\alpha}_{0}$ for $i=0$ as defined in (\ref{flux}), we also recover the result in that case.
\end{Proof}

\subsection{Numerical extension for non-fixed coefficients $(\gamma^\alpha)$}

Up to now, we were considering fixed coefficients $\gamma:=(\gamma^\alpha)_\alpha$ and the flux of the scheme at the junction point at time step $n \geq 0$ was
$$F_0 (\gamma) := \min \left\{ \displaystyle{\min_{\beta\leq N_I}} \ \dfrac{1}{\gamma^\beta} f_D^\beta(\rho^{\beta,n}_{-1}), \ \displaystyle{\min_{\lambda \geq N_I+1}} \ \dfrac{1}{\gamma^\lambda} f_S^\lambda(\rho^{\lambda,n}_{0}) \right\}.$$

In certain situations, we want to maximize the flux $F_0 (\gamma)$ for $\gamma$ belonging to an admissible set $\Gamma$. Indeed we can consider the set
$$A:= \argmax{\gamma \in \Gamma} F_0(\gamma).$$

In the case where this set is not a singleton, we can also use a priority rule to select a single element $\gamma^{*,n}$ of $A$.
This defines a map
$$\left( (\rho^{\beta,n}_{-1})_{\beta \leq N_I},(\rho^{\lambda,n}_{0})_{\lambda \geq N_I+1} \right) \mapsto \gamma^{*,n}.$$

At each time step $n \geq 0$ we can then choose this value $\gamma = \gamma^{*,n}$ in the numerical scheme (\ref{LWRscheme}), (\ref{flux}).

\section{Simulation} 
\label{simulation}

In this section, we present a numerical experiment. The main goal is to check if the numerical scheme (\ref{NS}),(\ref{CINS}) (or equivalently the scheme (\ref{LWRscheme}),(\ref{CI-LWR})) is able to illustrate the propagation of shock or rarefaction waves for densities on a junction.

\subsection{Settings}

We consider the case of a junction with $N_I=2=N_O$, that is two incoming roads denoted $\alpha=1$ and $2$ and two outgoing roads denoted $\alpha=3$ and $4$.

\begin{figure}[!ht]
\begin{center}
\includegraphics[scale=0.45]{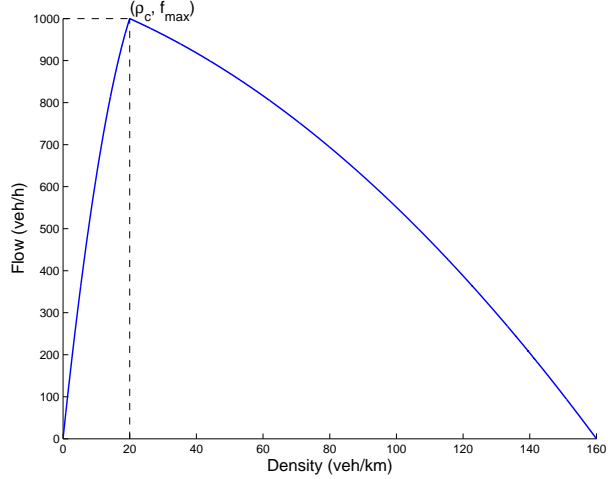}
\caption{Graph of the function $f$}
\label{fig::1}
\end{center}
\end{figure}

For the simulation, we consider that the flow functions are equal on each branch $f^{\alpha}=:f$ for any $\alpha=1,...,4$. Moreover the function $f$ is bi-parabolic (and only Lipschitz) as depicted on Figure \ref{fig::1}. It is defined as follows
$$f(\rho)=
\begin{cases}
\dfrac{(1-k) f_{max}}{\rho_c^2} \rho^2 + \dfrac{k f_{\max}}{\rho_c} \rho, \quad &\mbox{for} \quad \rho \leq \rho_c, \\
\\
\dfrac{(1-k) f_{max}}{(\rho_{max}-\rho_c)^2} \rho^2 + \dfrac{(k \rho_c + (k-2)\rho_{max}) f_{max}}{(\rho_{max}-\rho_c)^2} \rho - \dfrac{\rho_{max}(k \rho_c - \rho_{max}) f_{max}}{(\rho_{max}-\rho_c)^2}, \quad &\mbox{for} \quad \rho > \rho_c,
\end{cases}$$
with the jam density $\rho_c=~20 \ veh/km$, the maximal $\rho_{max}=~160 \ veh/km$, the maximal flow $f_{max}=~1000 \ veh/h$ and $k=1.5$.

The Hamiltonians $H^\alpha$ for $\alpha=1,...,4$ are defined in (\ref{Hfunctf}) according to the flow function $f$. See also Remark~\ref{weakerassumptions} on weaker assumptions than (A1) on the Hamiltonians. 
We also assume that the coefficients $(\gamma^\alpha)$ are all identical
$$\gamma^\alpha=\dfrac{1}{2} \quad \mbox{for any} \quad \alpha=1,...,4.$$

Notice that the computations are carried out for different $\Delta x$. In each case the time step $\Delta t$ is set to the maximal possible value satisfying the CFL condition (\ref{CFL'}). We consider branches of length $L=200 \ m$ and we have $\displaystyle N_b:=\left\lfloor \dfrac{L}{\Delta x} \right\rfloor$ points on each branch such that $i \in \{0,...,N_b \}$.

\subsection{Initial and boundary conditions}

\textbf{Initial conditions.} \ In traffic flow simulations it is classical to consider Riemann problems for the vehicles densities at the junction point. We not only consider a Riemann problem at the junction but we also choose the initial data with a second Riemann problem on the outgoing branch number 3 (see Table \ref{tab::1} where \textit{left} (resp. \textit{right}) stands for the left (resp. right) section of branch 3 according to this Riemann problem). We then consider initial conditions $(u^\alpha_0(x))_{\alpha=1,...,N}$ corresponding to the primitive of the densities depicted on Figure \ref{fig::2} (a). We also take the initial label at the junction point such that
$$u^{\alpha}_0(0)=:u_0(0)=0, \quad \mbox{for any} \quad \alpha.$$
We can check that the initial data $(u^\alpha_0(x))_{\alpha=1,...,N}$ satisfy (A0).

We are interested in the evolution of the densities. We stop to compute once we get a stationary final state as shown on Figure \ref{fig::2} (f). The values of densities and flows are summarized in Table \ref{tab::1}.

\begin{table}[h!]
\begin{center}
\begin{tabular}{ccccc}
 & \multicolumn{2}{c}{Initial state} & \multicolumn{2}{c}{Final state} \\
Branch & Density & Flow & Density & Flow \\
 & (veh/km) & (veh/h) & (veh/km) & (veh/h) \\ \hline
1 & 15 & \textit{844} & 90 & \textit{625} \\
2 & 15 & \textit{844} & 90 & \textit{625} \\
3 (left) & 30 & \textit{962} & 90 & \textit{625} \\
3 (right) & 90 & \textit{625} & 90 & \textit{625} \\
4 & 5 & \textit{344} & 10 & \textit{625}
\end{tabular}
\caption{Values of densities and flows for initial and final states on each branch}
\label{tab::1}
\end{center}
\end{table}

\textbf{Boundary conditions.} \ For any $i \leq N_b$ we use the numerical scheme (\ref{NS}) for computing $(U^{\alpha,n}_{i})$. Nevertheless at the last grid point $i=N_b$, we have
$$\dfrac{U^{\alpha ,n+1}_{N_b} -U^{\alpha ,n}_{N_b}}{\Delta t} + \max \left\{ H^{+}_{\alpha} (p^{\alpha ,n}_{N_b,-}), H^{-}_{\alpha} (p^{\alpha ,n}_{N_b,+}) \right\} =0, \quad \mbox{for}\quad  \alpha = 1,\ldots,N ,$$
where $p^{\alpha ,n}_{N_b,-}$ is defined in (\ref{defp}) and we set the boundary gradient as follows
$$p^{\alpha,n}_{N_b,+}=\begin{cases}
\dfrac{\rho^\alpha_0}{\gamma^{\alpha}}, \quad &\mbox{if} \quad \alpha \leq N_I, \\
p^{\alpha,n}_{N_b,-}, \quad &\mbox{if} \quad \alpha \geq N_I+1.
\end{cases}$$

These boundary conditions are motivated by our traffic application. Indeed while they are presented for the scheme (\ref{NS}) on $(U^{\alpha,n}_{i})$, the boundary conditions are easily translatable to the scheme (\ref{LWRscheme}) for the densities. For incoming roads, the flow that can enter the branch is given by the minimum between the supply of the first cell and the demand of the virtual previous cell which correspond to the value of $f$ evaluated for the initial density on the branch $\rho^\alpha_0$ (see Table \ref{tab::1}).
For outgoing roads, the flow that can exit the branch is given by the minimum between the demand of the last cell and the supply of the virtual next cell which is the  same than the supply of the last cell.

\subsection{Simulation results}

\textbf{Vehicles labels and trajectories.} \ Notice that here the computations are carried out for the discrete variables $(U^{\alpha,n}_i)$ while the densities $(\rho^{\alpha,n}_i)$ are computed in a post-treatment using (\ref{defrho}). It is also possible to compute directly the densities $(\rho^{\alpha,n}_i)$ according to the numerical scheme (\ref{LWRscheme}). Hereafter we consider $\Delta x= 5 m$ (that corresponds to the average size of a vehicle) and $\Delta t= 0.16 s$.

\begin{figure}[!ht]
 \begin{tabular}{cc}
 \includegraphics[scale=0.45]{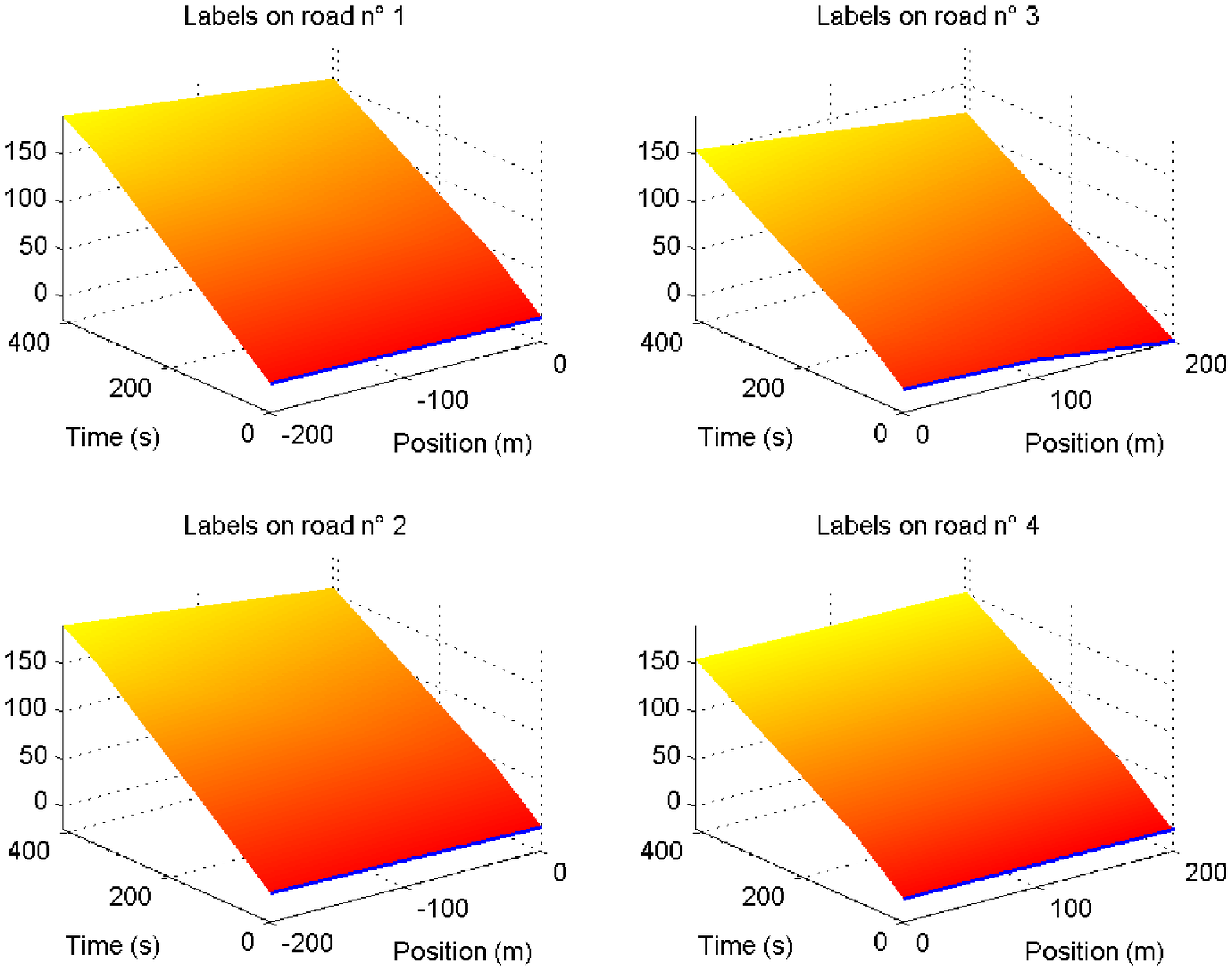} & \includegraphics[scale=0.45]{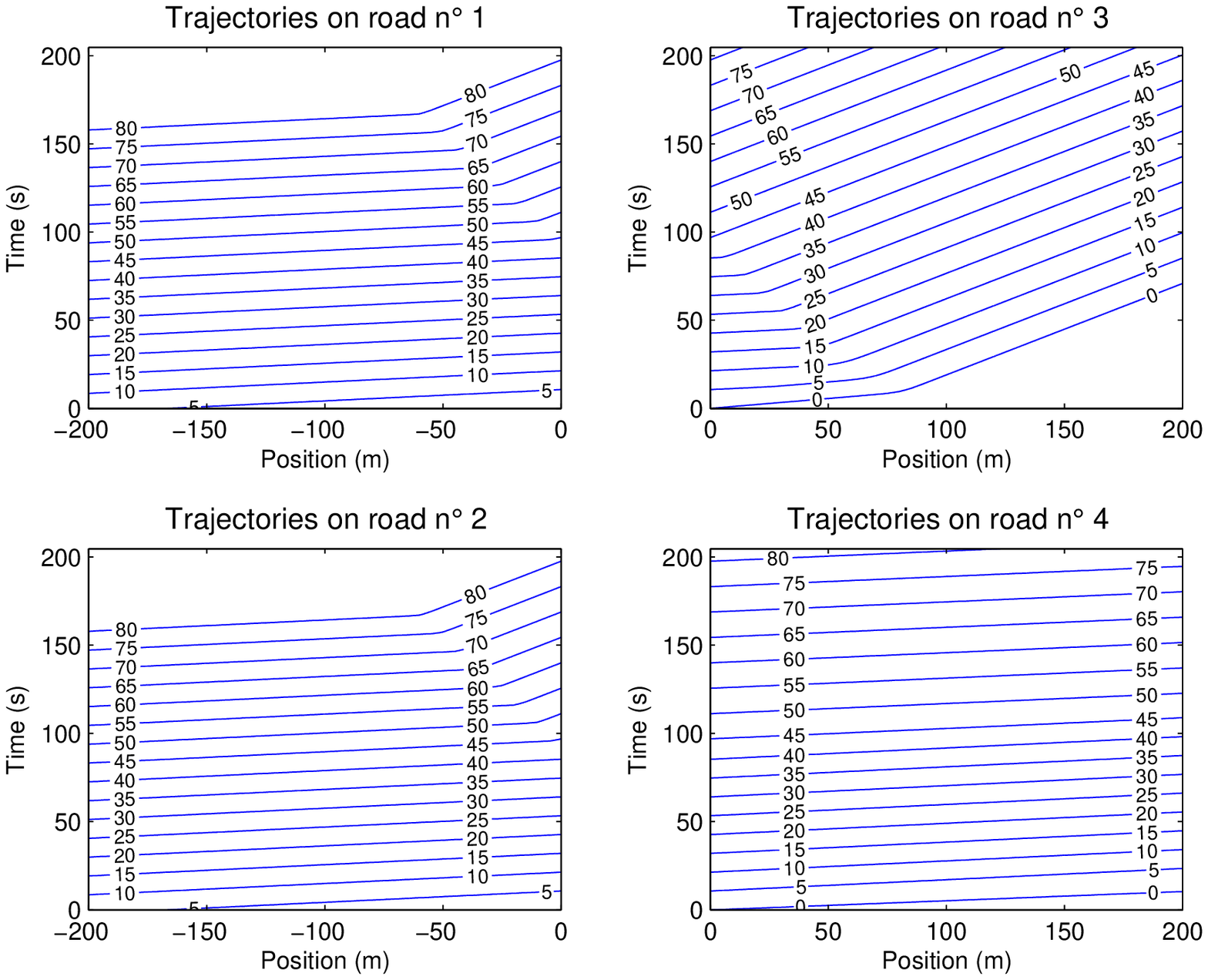} \\
 (a) Discrete labels $(U^{\alpha,n}_i)_{i,n}$ & (b) Trajectories of some vehicles
 \end{tabular}
 \caption{Numerical solution and vehicles trajectories}
 \label{fig::3}
 \end{figure}

The numerical solution $(U^{\alpha,n}_i)$ is depicted on Figure \ref{fig::3} (a). The vehicles trajectories are deduced by considering the iso-values of the labels surface $(U^{\alpha,n}_i)$ (see Figure \ref{fig::3} (b)). In this case, one can observe that the congestion (described in the next part) induces a break in the velocities of the vehicles when going through the shock waves. The same is true when passing through the junction.

\begin{figure}[!ht]
 \begin{center}
 \includegraphics[scale=0.75]{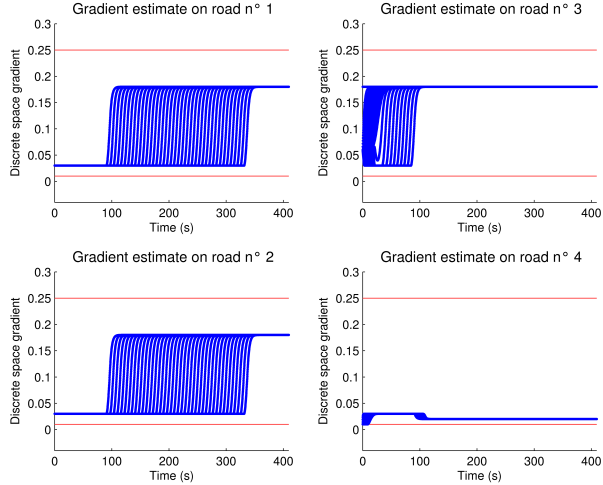}
 \caption{Bounds $\overline{p}^\alpha$ and $\underline{p}^\alpha$ on the gradient}
 \label{fig::4}
 \end{center}
 \end{figure}

We can also recover the gradient properties of Theorem~\ref{gradientestimate}. On Figure \ref{fig::4}, the gradients $\left( p^{\alpha,n}_{i,+} \right)$ are plotted as a function of time. We numerically check that the gradients stay between the bounds $\overline{p}^\alpha$ and $\underline{p}^\alpha$.\\

 \begin{figure}
 \begin{tabular}{ccc}
 \includegraphics[scale=0.40]{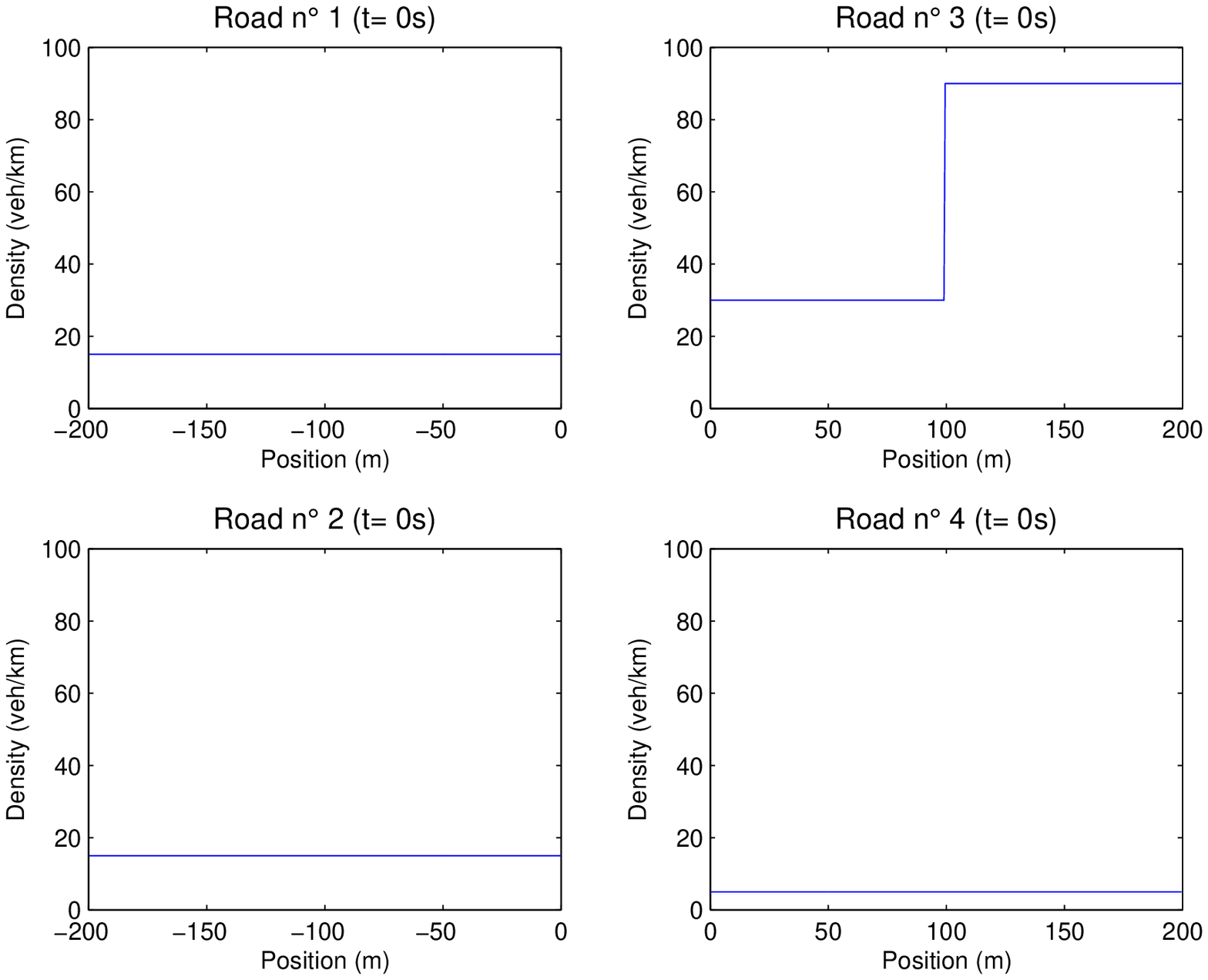} & & \includegraphics[scale=0.40]{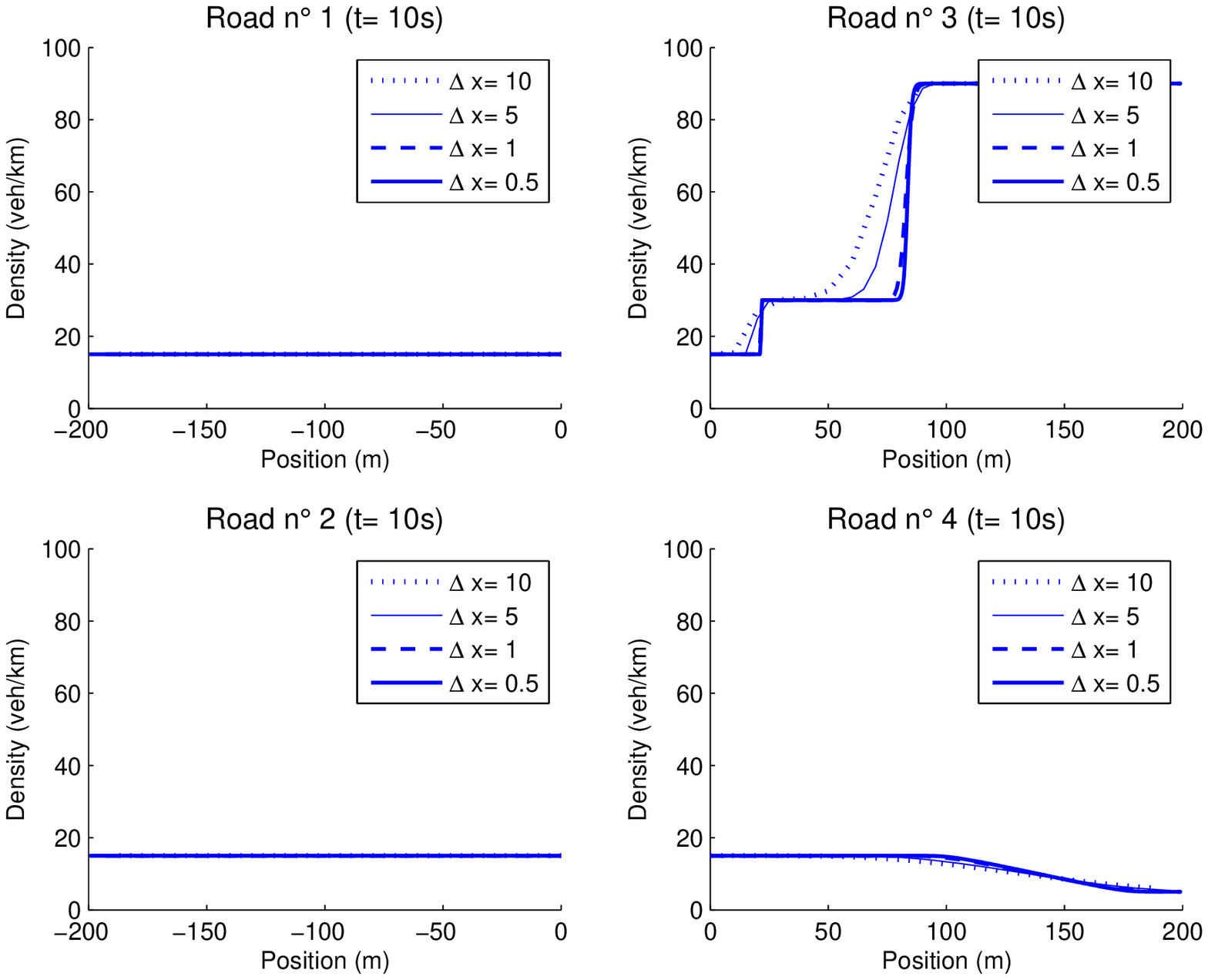} \\
 (a) Initial conditions for densities & & (b) Densities at $t=10$ s \\
  & & \\
  & & \\
 \includegraphics[scale=0.40]{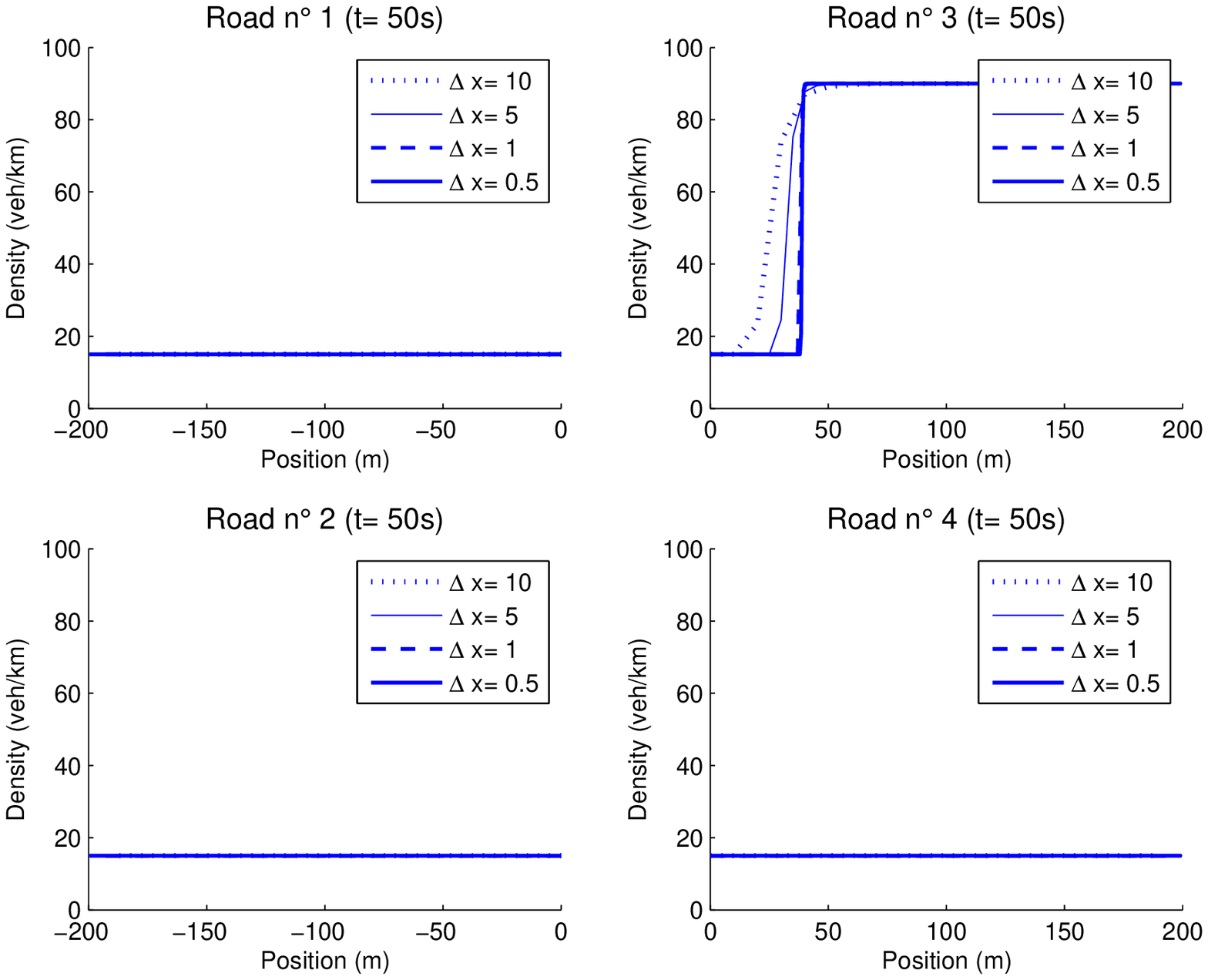} & & \includegraphics[scale=0.40]{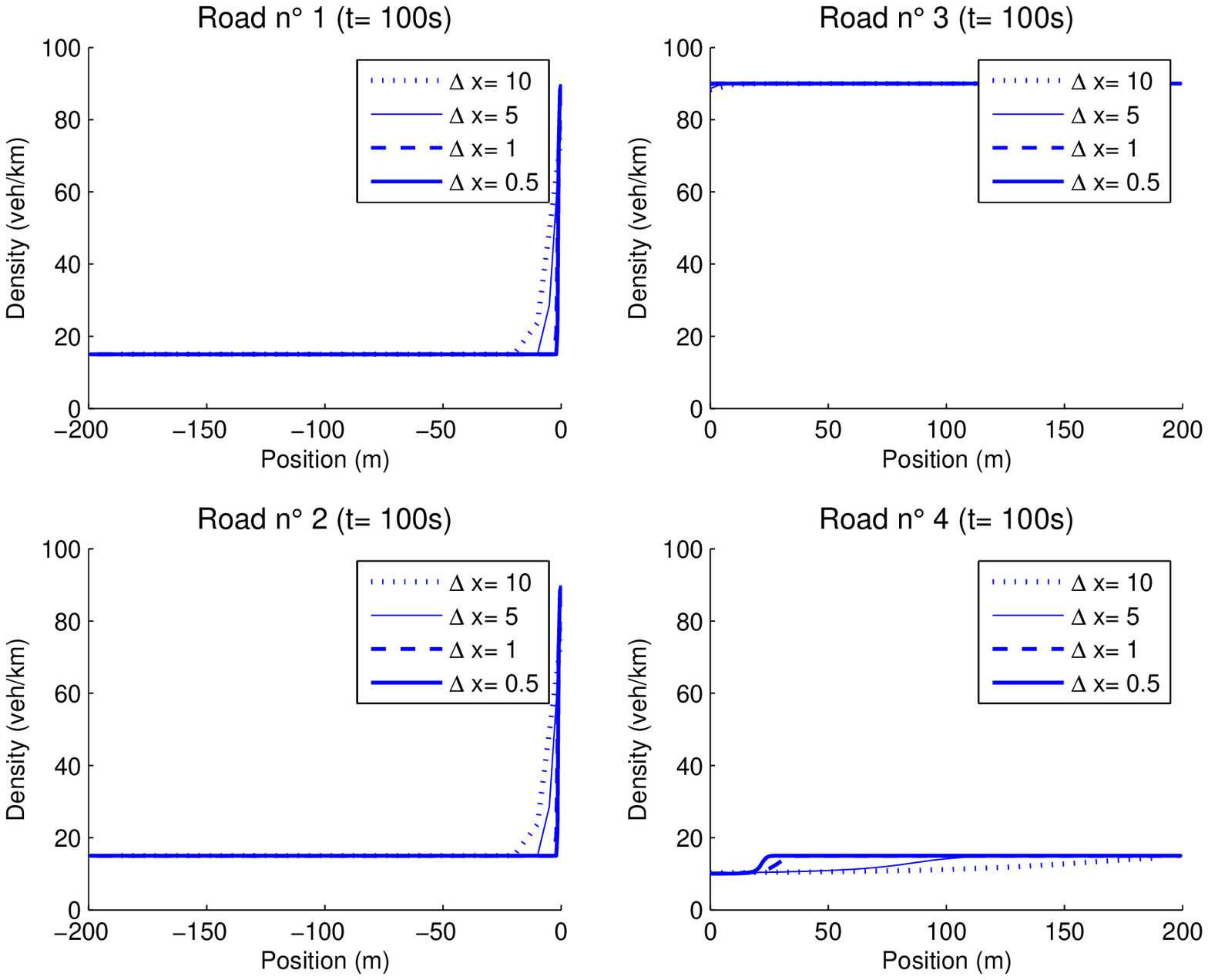} \\
 (c) Densities at $t=50$ s & & (d) Densities at $t=100$ s \\
  & & \\
  & & \\
 \includegraphics[scale=0.40]{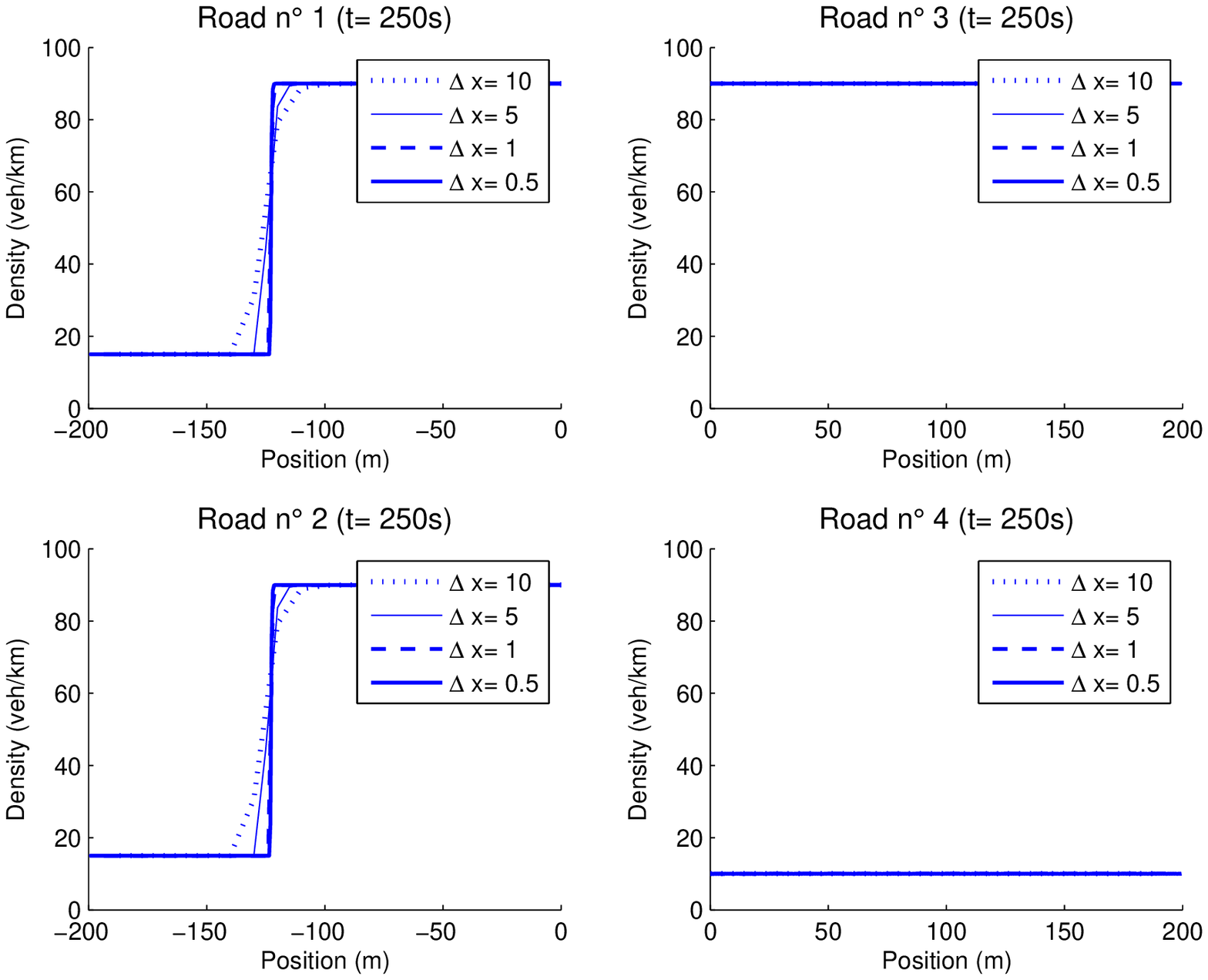} & & \includegraphics[scale=0.40]{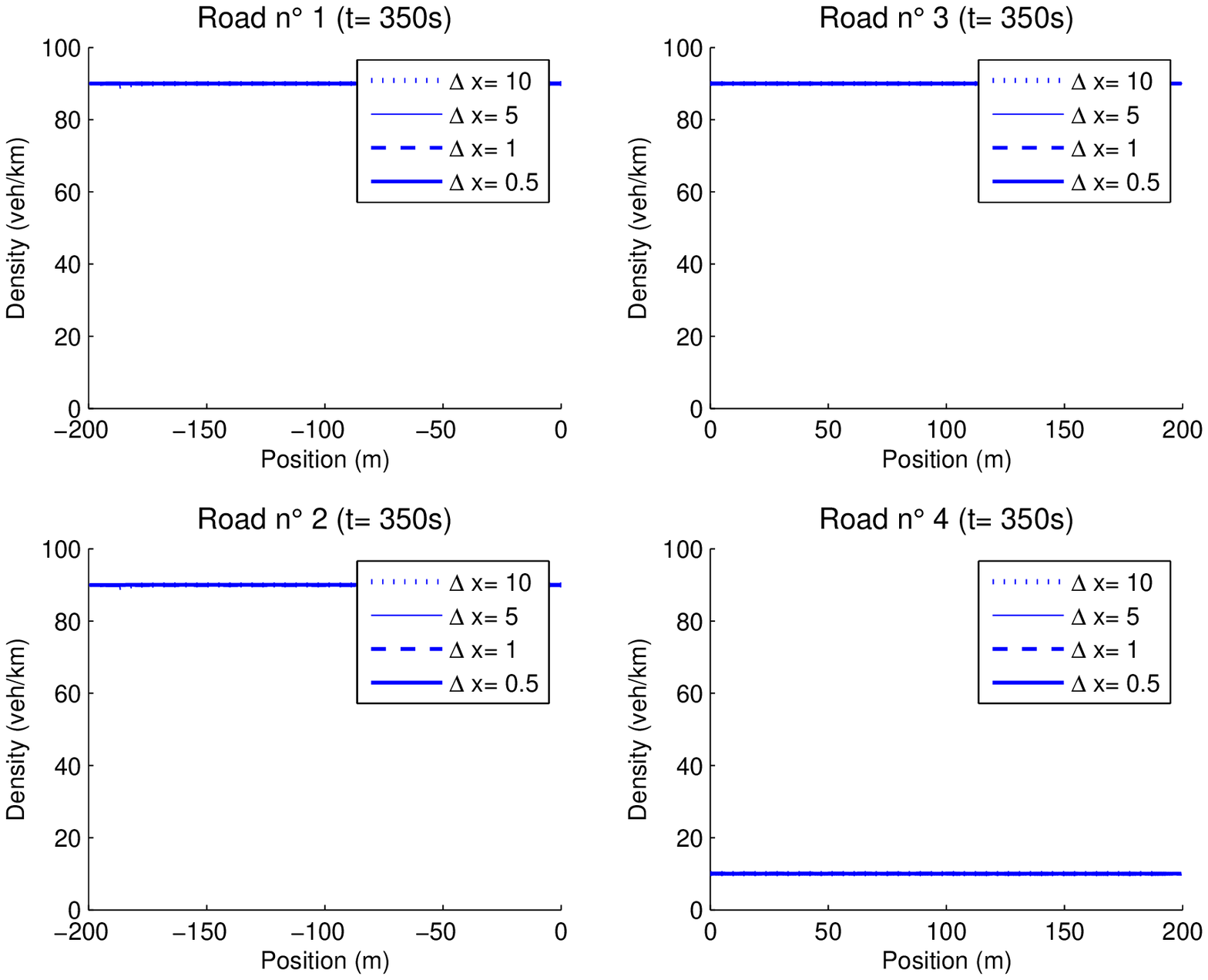} \\
 (e) Densities at $t=250$ s &  & (f) Densities at $t=350$ s \\
  & & 
 \end{tabular}
 \caption{Time evolution of vehicles densities for different $\Delta x$}
 \label{fig::2}
 \end{figure} 

\textbf{Propagation of waves.} \ We describe hereafter the shock and rarefaction waves that appear from the considered initial Riemann problems (see Figure \ref{fig::2}). At the initial state (Figure \ref{fig::2} (a)), the traffic situation on roads 1, 2 and 4 is fluid ($\rho^{\{1,2,4\}}_0 \leq \rho_c$) while the road 3 is congested ($\rho^3_0 \geq \rho_c$). Nevertheless the demands at the junction point are fully satisfied. As we can see on Figure \ref{fig::2} (b), there is the apparition of a rarefaction wave on road 4 and a shock wave on road 3, just downstream the junction point. At the same time, there is a shock wave propagating from the middle of the section on road 3 due to the initial Riemann problem there. This shock wave should propagate backward at the Rankine-Hugoniot speed $\tilde{v}_1=-6$ km/h. A while later (Figure \ref{fig::2} (c)), the rarefaction wave coming for the junction point and the shock wave coming from the middle of road 3 generate a new shock wave propagating backward at the speed of $\tilde{v}_2=-3$ km/h. The congestion spreads all over the branch 3 and reaches the junction point. At that moment (Figure \ref{fig::2} (d)), the supply on road 3 (immediately downstream the junction point) collapses. The demand for road 3 cannot be satisfied. Then it generates a congestion on both incoming roads. The shock wave continues to propagate backward in a similar way on roads 1 and 2 at speed $\tilde{v}_2$ (Figure \ref{fig::2} (e)). This congestion decreases also the possible passing flow from the incoming roads to the road 4. There is then a rarefaction wave that appears on road 4. However road 3 is still congested while the traffic situation on road 4 is fluid (Figure \ref{fig::2} (f)).\\

Figure \ref{fig::2} numerically illustrates the convergence of the numerical solution $(\rho^{\alpha,n}_i)$ when the grid size $(\Delta x, \Delta t)$ goes to zero. The rate of convergence is let to further research.

\section*{Aknowledgements}
The authors are grateful to C. Imbert for indications about the literature, M. Husti\'{c} for his suggestions to simplify certain parts of the proofs and L. Paszkowski for valuable comments on the presentation.

This work was partially supported by the ANR (Agence Nationale de
la Recherche) through HJnet project ANR-12-BS01-0008-01.


\end{document}